\documentclass{article}
\usepackage{geometry}
\usepackage{cite}
\usepackage{authblk}
\usepackage{amssymb,hyperref,amsmath,amsthm}
\usepackage{tikz}
\usetikzlibrary{decorations.pathreplacing,calc}
\usepackage{subcaption}
\pgfdeclarelayer{bg}    
\pgfsetlayers{bg,main}

\author[1, 2]{Ant\'onio Leit\~ao}
\author[1,3]{Nina Otter}

\affil[1]{ DataShape, Inria-Saclay,  France}
\affil[2]{Scuola Normale Superiore di Pisa, Italy}
\affil[3]{ Laboratoire de Mathématiques d’Orsay,  Universit\'e Paris-Saclay, Orsay, France}

\theoremstyle{definition}
\newtheorem{theorem}{Theorem}[section]

\newtheorem{definition}[theorem]{Definition}
\newtheorem{example}[theorem]{Example}

\title{Time-optimal persistent homology representatives for univariate time series} 

\newcommand{\define}[1]{{\bf \boldmath{#1}}}
\newcommand{\im}{\mathrm{im}}
\newcommand{\bd}{\partial}
\newcommand{\birth}{\mathrm{birth}}
\newcommand{\death}{\mathrm{death}}
\newcommand{\argmin}{\mathrm{argmin}}

\newcommand{\FF}{\mathbb{F}}
\date{}

\begin{document}

\maketitle

\begin{abstract}

Persistent homology (PH) is one of the main methods used in Topological Data Analysis. 
An active area of research in the field is the study of appropriate notions of PH representatives, which allow to interpret the meaning of the information provided by PH, making it an important problem in the application of PH, and in the study of its interpretability.
Computing optimal PH representatives is a problem that is known to be NP-hard, and one is therefore interested in developing context-specific optimality notions that are computable in practice. Here we introduce   time-optimal PH representatives for time-varying data, allowing one to extract representatives that are close in time in an appropriate sense.
We illustrate our methods on quasi-periodic synthetic time series, as well as time series arising from climate models, and we show that our methods provide optimal PH representatives that  are better suited for these types of problems than existing optimality notions, such as length-optimal PH representatives. 
\end{abstract}

\section{Introduction}

Topological Data Analysis (TDA) is a field that uses insights from topology --- the mathematical area that studies abstract shapes --- to develop representations of data that are computable and robust in an appropriate sense \cite{Carlsson2009TopologyAD}. 
Persistent homology (PH)  is, arguably, one of the most successful methods used in Topological Data Analysis, and it is being increasingly applied to a variety of data analysis problems. We refer the reader to the DONUT database \cite{DONUT} for a vast collection or real-world applications of PH. In persistent homology, one takes as a point of departure a data set, such as a point cloud, a time series, a network, or a digital image, and associates to it a $1$-parameter family of topological spaces, in which for each parameter value $r$ one may think of the corresponding space as being an approximation, truncation or thickening of the original data set; for instance, if the parameter captures distance between points, the space at parameter value $r$ identifies all points at distance smaller or equal than $r$. The output of persistent homology is then  a summary, called ``persistence barcode'' or ``persistence diagram'', of the number of topological features such as connected components, holes or tunnels, voids, present in a data set, as well as how long each feature spans (``persists'') across the possible parameter values. 

Often, in application, one is  interested in interpreting the meaning of such summaries, by determining which data points correspond to each individual summary. More precisely, given a persistence barcode, one asks for a choice of persistent vector basis, called ``PH representatives'', that is meaningful for the application at hand. 
Ideally, one would want to compute  representatives that are optimal in an appropriate sense. This problem has been shown to be NP hard \cite{CF11}, and a considerable amount of work is being devoted to developing algorithms to find approximations of representatives that satisfy some minimality condition of interest in a specific application context.

\begin{figure}
\[
\includegraphics[scale=0.2]{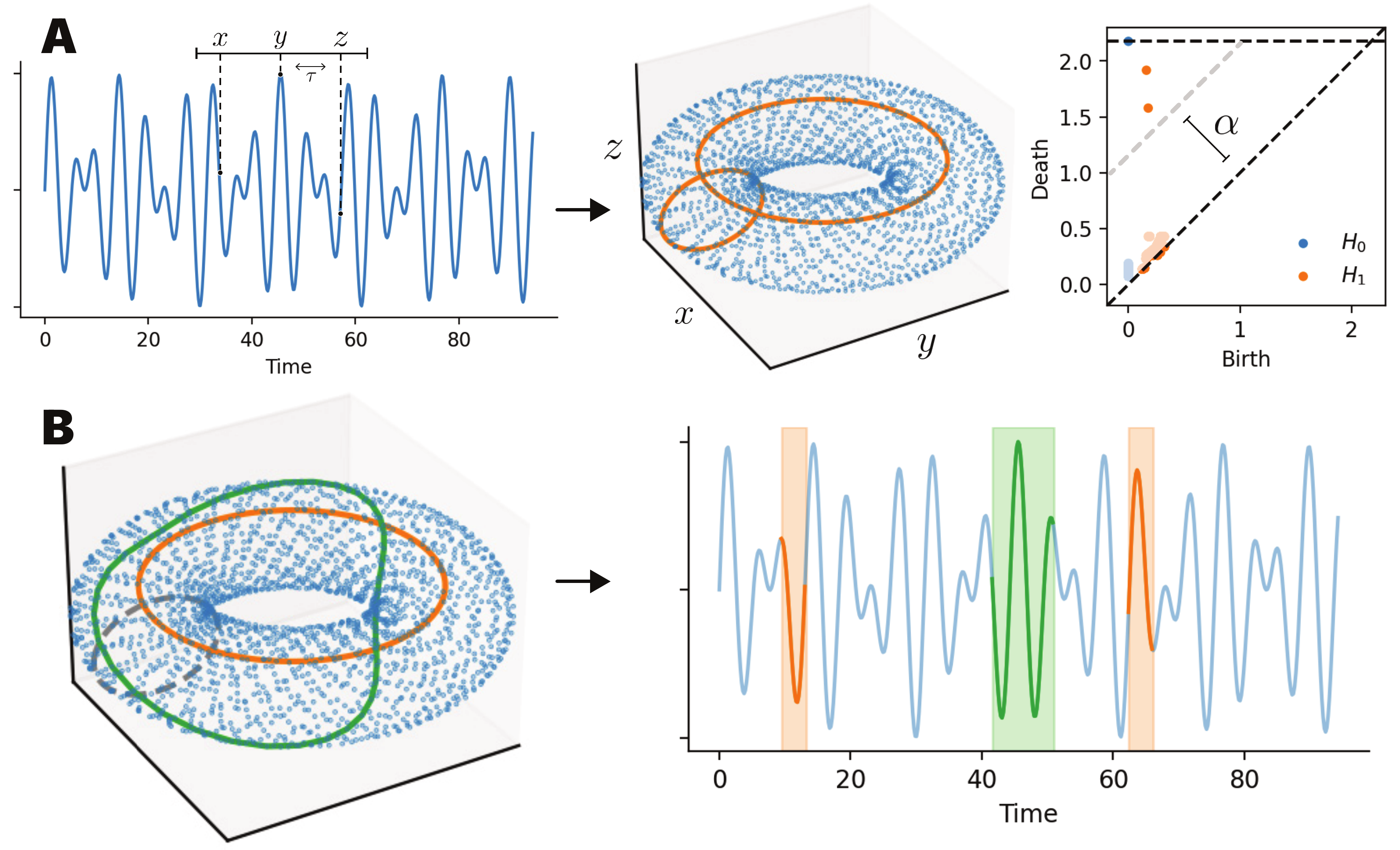}
\]
 \caption{\textbf{Pipeline for computing time-optimal PH cycle representatives for univariate time series}. \textbf{A)} One starts by embedding a univariate time series into Euclidean space. For quasi-periodic times series, one can obtain an embedding together with a lower bound ($\alpha$) on  persistence that allows to distinguish between components, holes, voids,  that one may consider as ``significant'', and   features that may be due to noise.  
 \textbf{B)} Given the significant PH features, we attempt to find representatives that are optimal with respect to their cohesion in time. Given representatives of the same PH feature (e.g., green and orange cycles on the left), the ones that correspond to a continuous trajectory (green) in the original signal are preferred over the ones that are discontinuous (orange).} \label{fig:pipeline}
\end{figure}

In the present work we introduce algorithms to compute PH representatives for time-dependent data. To be best of our knowledge, this is the first time that such representatives have been studied. 
In particular, we propose two different notions: (i) vertex time-optimal PH representatives and (ii) simplex time-optimal PH representatives. We illustrate these notions on synthetic quasi-periodic time series and on  time series resulting from delayed oscillator models of the El Nin\~o Southern Oscillation (ENSO). 
We illustrate the pipeline for the computation of time-optimal representatives for univariate time series
in Figure \ref{fig:pipeline}.

\section{Related work}

To the best of our knowledge, our work is the first attempt at studying PH representatives for time-varying data that explicitly take the time information into account. 
Our work inscribes itself in a line of research that tries to develop notions of optimal PH representatives  suitable for applications, such as the length-optimal or volume-optimal PH representatives \cite{DHM19,obayashi2018volume}. A benchmarking of several existing approaches has been performed in \cite{lu2021minimal}.
A speed-up exploiting persistent cohomology computations has been proposed in \cite{CV23}.

The development of  methods to study time-varying data is a very active area of research in TDA. Vineyards provide $1$-parameter families of persistence barcodes for time-varying data \cite{CSEM06}; generalisations of these to multi-parameter families have recently been proposed \cite{H22b, H22a}. An algebraic framework for the study of time-varying persistence modules has been introduced in \cite{T23}.
Univariate time series  have been studied in TDA, among others, in \cite{kenneth2009human,perea2015sw1pers,RMC23,PH15 }.
Several topological methods have been developed to study bifurcation diagrams of dynamical systems, including \cite{GMK22}.

\section{Background}
We first  give a brief overview of basic notions from  persistent homology in Section \ref{SS:PH repr}; we then discuss univariate time series in Section \ref{SS:time-var data}, and we give an overview on existing optimization approaches for PH cycle representatives in Section \ref{SS:optim PH cycles}.

\subsection{Homology and persistent homology representatives}\label{SS:PH repr}

To simplify exposition, here we introduce homology and persistent homology for coefficients over the field with two elements $\mathbb{F}_2$. We refer the reader to Appendix \ref{A:PH} for a discussion of what changes for arbitrary coefficient fields.

\subsubsection{Homology}
 Let $K$ be a simplicial complex,  and denote by $S_p(K)$ its set of $p$-simplices. 
We denote by  $C_p(K)$  the vector space generated by the $\mathbb{F}_2$-linear combinations of $p$-simplices.

The elements of $C_p(K)$ are called \define{$p$-chains}.
We consider the boundary operator 
\begin{alignat}{2}
\bd_p: C_p(K) &\notag\longrightarrow C_{p-1}(K)\\
\sigma &\notag \mapsto \sum_{\tau \subset \sigma, \tau\in S_{p-1}(K)}  \tau
\end{alignat}

and call the elements of the kernel of $\bd_p$ \define{$p$-cycles},   
and the elements of the image of $\bd_{p+1}$ \define{$p$-boundaries}. One can show that $\bd_{p+1}\circ \bd_p=0$ for all $p$. Intuitively, this is due to the fact that the boundary of a boundary is empty.
The quotient vector space $H_p(K) = \frac{\ker(\bd_p)}{\im (\bd_{p+1})}$ is called \define{$p$th simplicial homology} of $K$ (with coefficients in $\FF_2$).

\begin{definition}
For a given simplicial homology vector space $H_p(K) = \frac{\ker(\bd_p)}{\im (\bd_{p+1})}$ and $h\in H_p(K)$, we call $c\in \ker(\bd_p)$ a \define{cycle representative} of $h$ if $[c]:=c+\im (\bd_{p+1}=h$. Two cycles $c,c'\in \ker(\bd_p)$ are \define{homologous} if $[c]=[c']$.
A \define{cycle basis} for $H_p(K)$ is a set of  $p$-cycles $c_1,\dots , c_m$ such that $[c_i]\ne [c_j]$ for $i\ne j$ and $[c_1],\dots , [c_m]$ are a basis for $H_p(K)$.

\end{definition}

\subsubsection{Persistent homology}\label{SS:PH}

We now consider a finite sequence of nested simplicial complexes:
\[
 \varnothing \subseteq K_{t_0} \subseteq K_{t_1} \subseteq \dots \subseteq K_{t_{n-1}} \subseteq K_{t_n} =: K \, .
\]
We call $\{K_{t_i}\}_{i=0}^n$ a \define{filtration} of $K$.
We obtain injective linear maps
$\iota_{i,j}\colon C_p(K_{t_i})\hookrightarrow C_p(K_{t_{j}})$ for all $0\leq i<j\leq n$, induced by the inclusions of simplicial complexes.
We can thus identify each $C_p(K_{t_i})$ with a vector subspace of $C_p(K_{t_j})$ for any $i<j$.

\begin{definition}
Given $c\in C_p(K)$, we define its \define{birth} to be 
\[
\birth(c)=\min\{i\mid c\in C_p(K_i)\}
\]
and its \define{death} to be 
\[
\death(c)=\min\left \{i \mid c\in \im\left(\bd_{p+1}\big|_{(C_{p+1}(K_i))}\right)\right\}
\]
where we use the convention  $\min \emptyset =\infty$.
\end{definition}

Similarly as for  chain vector spaces, the inclusions of simplicial complexes induce (not necessarily injective) linear maps between homology vector spaces $\phi_{i,j}\colon H_p(K_{t_i})\to H_p(K_{t_j})$ for all $0\leq i\leq j \leq n$. One gives the following definition:
\begin{definition}
  The  \define{$p$th persistent homology} $H_p(K)$ of a filtered simplicial complex 
  $\{K_{t_i}\}_{i=0}^n$ with $K_{t_n}=K$ is the tuple 
  $(\{H_p(K_{t_i})\}_{i=0}^n,\{\phi_{i,j}\}_{i<j})$.
\end{definition}

\begin{definition}
A \define{persistent homology cycle basis} of $H_p(K)$ is a set of $p$-cycles $c_1,\dots , c_m\in C_p(K)$ such that $\birth(c_j)\ne \death(c_j)$ for all $j=1,\dots , m$ and for each filtration value $t$ we have that the collection of cycles $c_j$ with $\birth(c_j)\leq t\leq \death(c_j)$ form a cycle basis for $H_p(K_{t})$. We say that a $p$-cycle $c\in C_p(K)$ is a \define{persistent homology cycle representative} for $H_p(K)$ if it is an element of a persistent homology cycle basis of $H_p(K)$.
\end{definition}

\begin{definition}
Let  $c_1,\dots , c_m\in C_p(K)$  be a persistent homology cycle basis of $H_p(K)$. We call the multiset
\[
PD_p(K):=\Big\{\big(\birth(c_j),\death(c_j)\big) \mid  j=1,\dots , m \Big\}\]
the \define{persistence diagram of $H_p(K)$, or the $p$th persistence diagram of $K$.}
We call the number $\death(c_j)-\birth(c_j)$ the \define{persistence} of $c_j$.
\end{definition}

It is a fundamental result in persistent homology that a persistent homology cycle basis exists for any persistent homology tuple satisfying appropriate finiteness conditions, and that the persistence diagram does not depend on the choice of PH cycle basis \cite{CZ05}.

\begin{example}\label{E:ex cycle repr}
Consider the  filtration of simplicial complexes $K_1\subset K_2:=K$ illustrated in Figure \ref{F:example disjoint support PH repr}(ii). We can think of this filtration as the  triangulation of a bent cylinder, with its two sides on the bottom, which we then scan from bottom to top, see Figure \ref{F:example disjoint support PH repr}(i). 
We can compute the $1$st simplicial homology of $K_2$, and we obtain that $H_1(K_2)\cong \FF_2$.
Some possible choices of $1$-cycle representatives include for instance $<ab>+<bc>+<ac>$,  $<a'b'>+<b'c'>+<a'b'>$ or $<ab'>+<b'b>+<bc>+<ca>$, where we denote by $<x_0 x_1>$ the vector  corresponding to the  $1$-simplex with vertices $x_0, x_1$. We depict these three different choices in Figure \ref{F:example disjoint support PH repr}(iii). In particular, we note that cycle representatives are in general not unique.

\end{example}
\begin{example}
We now use the same filtered complex from Example \ref{E:ex cycle repr}  to  illustrates PH cycle representatives, as well as some subtleties in their interpretation. 
The persistence diagram of $H_1(K)$ contains two points not on the diagonal: the point $(1,2)$ and the point $(1,\infty)$. A possible choice of  PH cycle basis of $H_1(K)$ is given by ${<a'b'>+<b'c'>+<a'c'>}$ and $ <a'b'>+<b'c'>+<a'c'>+<ab>+<bc>+<ac>$. We illustrate the simplices corresponding to the two cycle representatives in Figure \ref{F:example disjoint support PH repr} in orange and in cyan, respectively. 
We note that the second PH cycle representative is an example of  representative with disjoint ``support''; namely, if we consider the subcomplex corresponding to  the $1$-simplices with non-zero coefficients, then it is disconnected. This defies the usual interpretation of PH $1$-cycle representatives as correponding to single ``tunnels'' or ``loops'' in the data.  
Interestingly, such PH cycle representatives do not seem to appear often in practice, see \cite[Section 6.6.1]{lu2021minimal}.

\begin{figure}[h!]
\centering
(i)
\includegraphics[scale=0.2]{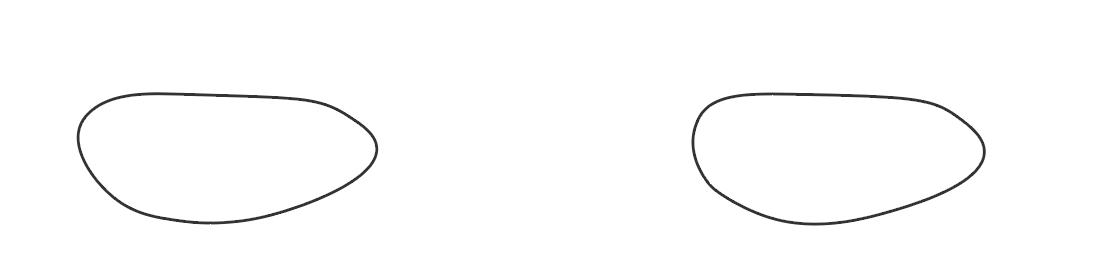}
$\subset \quad$
\includegraphics[scale=0.2]{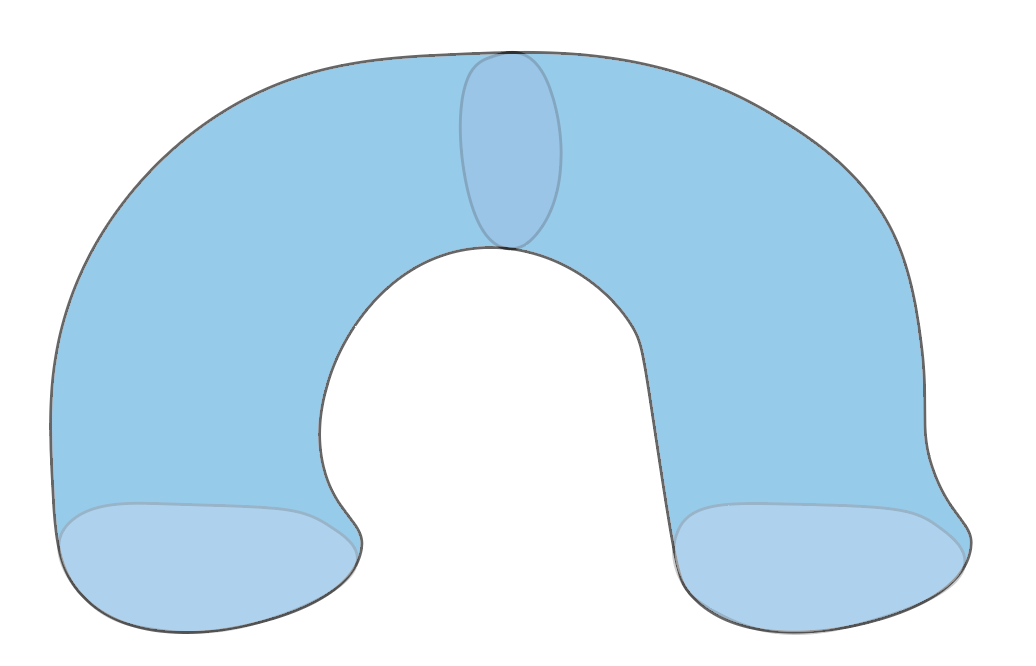}\\
\vspace{0.5cm}

(ii)\begin{tikzpicture}
\node at (1,-2.5) {$K_1$};
\node at (5.5,-2.5) {$K_2$};
\node at (0,0){$
\begin{tikzpicture}
\node (a) at (0,0) {$a$};
\node (b)  at (0.5,1) {$b$};
\node (c) at (1,0) {$c$};
\draw[-] (a)--(b)--(c)--(a);
\node (a') at (2,0) {$a'$};
\node (b') at (2.5,1) {$b'$};
\node (c') at (3,0) {$c'$};
\draw[-] (a')--(b')--(c')--(a');
\end{tikzpicture}$};
\node at (2.5,0)
{$\subset$};
\node at (5,0){$
 \begin{tikzpicture}
\draw[-,fill=yellow!20,draw opacity=0] (0,-3)--(2,-2)--(2,-3)--(0,-3)-- cycle;
\draw[-,fill=yellow!20,draw opacity=0] (0,-3)--(0,-2)--(2,-2)--(0,-3)-- cycle;
\draw[-,fill=yellow!20,draw opacity=0] (0,-2)--(2,-1)--(2,-2)--(0,-2)-- cycle;
\draw[-,fill=yellow!20,draw opacity=0] (0,-2)--(0,-1)--(2,-1)--(0,-2)-- cycle;
\draw[-,fill=yellow!20,draw opacity=0] (0,-1)--(2,0)--(2,-1)--(0,-1)-- cycle;
\draw[-,fill=yellow!20,draw opacity=0] (0,-1)--(0,0)--(2,0)--(0,-1)-- cycle;
\node (a) at (0,0) {$a$};
\node (c) at (0,-1) {$c$};
\node (b) at (0,-2) {$b$};
\node (a2) at (0,-3) {$a$};
\node (a') at (2,0) {$a'$};
\node (c') at (2,-1) {$c'$};
\node (b') at (2,-2) {$b'$};
\node (a'2) at (2,-3) {$a'$};
 \draw[-] (a)--(c)--(b)--(a2);
 \draw[-] (a')--(c')--(b')--(a'2);
\draw[-] (a2)--(a'2);
\draw[-] (a2)--(b');
\draw[-] (b)--(b');
\draw[-] (b)--(c');
\draw[-] (c)--(c');
\draw[-] (c)--(a');
\draw[-] (a)--(a');
\draw[-, thick, orange, dashed] (a')--(c')--(b')--(a'2);

\end{tikzpicture}
$};
\end{tikzpicture}\\
\vspace{0.5cm}
(iii)\qquad
\begin{tikzpicture}
\draw[-,fill=yellow!20,draw opacity=0] (0,-3)--(2,-2)--(2,-3)--(0,-3)-- cycle;
\draw[-,fill=yellow!20,draw opacity=0] (0,-3)--(0,-2)--(2,-2)--(0,-3)-- cycle;
\draw[-,fill=yellow!20,draw opacity=0] (0,-2)--(2,-1)--(2,-2)--(0,-2)-- cycle;
\draw[-,fill=yellow!20,draw opacity=0] (0,-2)--(0,-1)--(2,-1)--(0,-2)-- cycle;
\draw[-,fill=yellow!20,draw opacity=0] (0,-1)--(2,0)--(2,-1)--(0,-1)-- cycle;
\draw[-,fill=yellow!20,draw opacity=0] (0,-1)--(0,0)--(2,0)--(0,-1)-- cycle;
\node (a) at (0,0) {$a$};
\node (c) at (0,-1) {$c$};
\node (b) at (0,-2) {$b$};
\node (a2) at (0,-3) {$a$};
\node (a') at (2,0) {$a'$};
\node (c') at (2,-1) {$c'$};
\node (b') at (2,-2) {$b'$};
\node (a'2) at (2,-3) {$a'$};
 \draw[-] (a)--(c)--(b)--(a2);
 \draw[-] (a')--(c')--(b')--(a'2);
\draw[-] (a2)--(a'2);
\draw[-] (a2)--(b');
\draw[-] (b)--(b');
\draw[-] (b)--(c');
\draw[-] (c)--(c');
\draw[-] (c)--(a');
\draw[-] (a)--(a');
\draw[-, ultra thick, orange] (a')--(c')--(b')--(a'2);
\draw[-, ultra thick, cyan] (a)--(c)--(b)--(a2);
\end{tikzpicture}
\qquad\qquad\qquad
\begin{tikzpicture}
\draw[-,fill=yellow!20,draw opacity=0] (0,-3)--(2,-2)--(2,-3)--(0,-3)-- cycle;
\draw[-,fill=yellow!20,draw opacity=0] (0,-3)--(0,-2)--(2,-2)--(0,-3)-- cycle;
\draw[-,fill=yellow!20,draw opacity=0] (0,-2)--(2,-1)--(2,-2)--(0,-2)-- cycle;
\draw[-,fill=yellow!20,draw opacity=0] (0,-2)--(0,-1)--(2,-1)--(0,-2)-- cycle;
\draw[-,fill=yellow!20,draw opacity=0] (0,-1)--(2,0)--(2,-1)--(0,-1)-- cycle;
\draw[-,fill=yellow!20,draw opacity=0] (0,-1)--(0,0)--(2,0)--(0,-1)-- cycle;
\node (a) at (0,0) {$a$};
\node (c) at (0,-1) {$c$};
\node (b) at (0,-2) {$b$};
\node (a2) at (0,-3) {$a$};
\node (a') at (2,0) {$a'$};
\node (c') at (2,-1) {$c'$};
\node (b') at (2,-2) {$b'$};
\node (a'2) at (2,-3) {$a'$};
 \draw[-] (a)--(c)--(b)--(a2);
 \draw[-] (a')--(c')--(b')--(a'2);
\draw[-] (a2)--(a'2);
\draw[-] (a2)--(b');
\draw[-] (b)--(b');
\draw[-] (b)--(c');
\draw[-] (c)--(c');
\draw[-] (c)--(a');
\draw[-] (a)--(a');
\draw[-, ultra thick, purple] (a2)--(b')--(b)--(c)--(a);
\end{tikzpicture}\\
\vspace{0.5cm}
(iv)\qquad
 \begin{tikzpicture}
\draw[-,fill=yellow!20,draw opacity=0] (0,-3)--(2,-2)--(2,-3)--(0,-3)-- cycle;
\draw[-,fill=yellow!20,draw opacity=0] (0,-3)--(0,-2)--(2,-2)--(0,-3)-- cycle;
\draw[-,fill=yellow!20,draw opacity=0] (0,-2)--(2,-1)--(2,-2)--(0,-2)-- cycle;
\draw[-,fill=yellow!20,draw opacity=0] (0,-2)--(0,-1)--(2,-1)--(0,-2)-- cycle;
\draw[-,fill=yellow!20,draw opacity=0] (0,-1)--(2,0)--(2,-1)--(0,-1)-- cycle;
\draw[-,fill=yellow!20,draw opacity=0] (0,-1)--(0,0)--(2,0)--(0,-1)-- cycle;
\node (a) at (0,0) {$a$};
\node (c) at (0,-1) {$c$};
\node (b) at (0,-2) {$b$};
\node (a2) at (0,-3) {$a$};
\node (a') at (2,0) {$a'$};
\node (c') at (2,-1) {$c'$};
\node (b') at (2,-2) {$b'$};
\node (a'2) at (2,-3) {$a'$};
 \draw[-] (a)--(c)--(b)--(a2);
 \draw[-] (a')--(c')--(b')--(a'2);
\draw[-] (a2)--(a'2);
\draw[-] (a2)--(b');
\draw[-] (b)--(b');
\draw[-] (b)--(c');
\draw[-] (c)--(c');
\draw[-] (c)--(a');
\draw[-] (a)--(a');
\draw[-, ultra thick, orange] (a')--(c')--(b')--(a'2);

\end{tikzpicture} \qquad\qquad\qquad
 \begin{tikzpicture}
\draw[-,fill=yellow!20,draw opacity=0] (0,-3)--(2,-2)--(2,-3)--(0,-3)-- cycle;
\draw[-,fill=yellow!20,draw opacity=0] (0,-3)--(0,-2)--(2,-2)--(0,-3)-- cycle;
\draw[-,fill=yellow!20,draw opacity=0] (0,-2)--(2,-1)--(2,-2)--(0,-2)-- cycle;
\draw[-,fill=yellow!20,draw opacity=0] (0,-2)--(0,-1)--(2,-1)--(0,-2)-- cycle;
\draw[-,fill=yellow!20,draw opacity=0] (0,-1)--(2,0)--(2,-1)--(0,-1)-- cycle;
\draw[-,fill=yellow!20,draw opacity=0] (0,-1)--(0,0)--(2,0)--(0,-1)-- cycle;
\node (a) at (0,0) {$a$};
\node (c) at (0,-1) {$c$};
\node (b) at (0,-2) {$b$};
\node (a2) at (0,-3) {$a$};
\node (a') at (2,0) {$a'$};
\node (c') at (2,-1) {$c'$};
\node (b') at (2,-2) {$b'$};
\node (a'2) at (2,-3) {$a'$};
 \draw[-] (a)--(c)--(b)--(a2);
 \draw[-] (a')--(c')--(b')--(a'2);
\draw[-] (a2)--(a'2);
\draw[-] (a2)--(b');
\draw[-] (b)--(b');
\draw[-] (b)--(c');
\draw[-] (c)--(c');
\draw[-] (c)--(a');
\draw[-] (a)--(a');
\draw[-, ultra thick , cyan] (a')--(c')--(b')--(a'2);
\draw[-, ultra thick, cyan] (a)--(c)--(b)--(a2);

\end{tikzpicture}

\caption{(i) A nested sequence of subspaces of $\mathbb{R}^3$. (ii) A filtered simplicial complex that is obtained as  a triangulation of the spaces in (i). 
(iii) Three Possible choices of $1$-cycle representatives for $H_1(K_2)$. (iv) PH $1$-cycle representative basis for $H_1(K)$: we depict the two PH cycle representatives in orange and cyan. This is an example of a filtered simplicial complex for which one of the PH $1$-cycle representatives consists of a linear combination of two disjoint closed curves.  }\label{F:example disjoint support PH repr}
\end{figure}
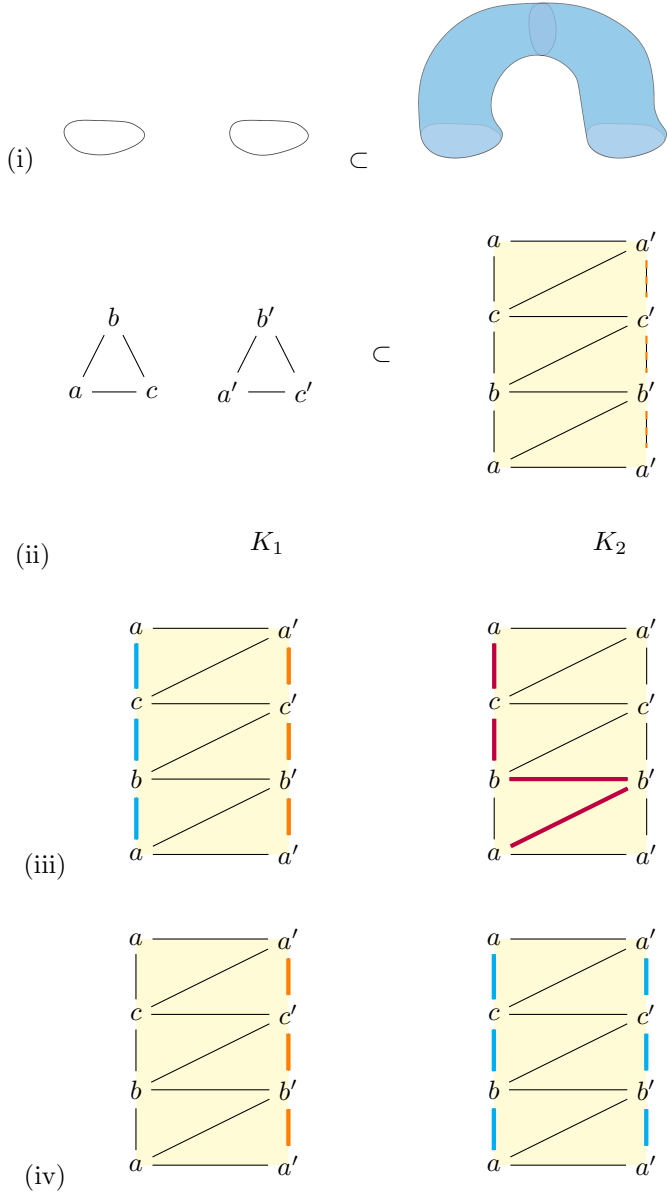
\end{example}

 \subsection{Univariate time series}\label{SS:time-var data}

\begin{definition}
A \define{univariate time series} is a function
 $f\colon \mathbb{R}\to \mathbb{R}$. 
\end{definition}

There are mainly two approaches for associating filtrations of simplicial complexes to  time series in TDA: (i) one computes sublevel-set filtrations of the time series or a transform thereof (such as a Discrete Fourier transform); or (ii) one embeds the time series in Euclidean space and associates filtrations of simplicial complexes to the resulting point cloud, such as, e.g., Vietoris-Rips complexes. See the review \cite{RC21} for details about these two approaches. 
Here we follow the second approach, and we provide details about the computation of the embedding parameters in  Appendix  \ref{A:optimal par}.

 We note that while in the current work we focus on univariate time series, our methods can be applied to a broader class of time-varying data, see  the discussion in Section \ref{S: conclusion}.

\subsection{Optimizing (persistent) homology cycle representatives}\label{SS:optim PH cycles}

In practice, for applications, we are interested in finding  cycle representatives that are informative for the problem at hand. In particular, one often seeks to find cycles that minimise some criteria (e.g., the number of simplices contained in the cycle). 
Roughly, the existing approaches to computing optimal PH $p$-cycle representatives can be divided into those that minimise a loss function defined on $p$-chains (also called  ``edge''-loss methods, in analogy with the $p=1$ case) and those that instead minimize a loss function defined on $p+1$-chains (also called ``triangle''-loss methods). We do not consider triangle-loss methods in the current work, and we point the reader to \cite{obayashi2018volume} for details. In the follow we review the basic set up of edge-loss methods.

\subsubsection{Edge-loss methods}\label{SSS:edge-loss}

Given an initial cycle $c_0 \in C_p(K)$ the problem for homology cycle representatives focuses on finding a homologous cycle that minimises a given loss function $\ell\colon C_{p}(K)\to \mathbb{R}$.
Since adding any boundary $w \in \im (\bd_{p+1})$ to $c_0$ results in a homologous cycle, we have the following  problem formulation:
\begin{align*}
\text{min}   & \quad \ell(c) \\
\text{subject to} & \quad c = c_0 + \bd_{p+1}(w) \\
                  & \quad w \in C_{p+1}(K)
\end{align*}

\noindent
In practice, obtaining the initial cycle is not a problem since most programs output a  cycle basis.

We can modify the previous problem to find a cycle that has the same persistence as the homology class it represents. Given a point $(b,d)$ in $PD_p(K)$ we look for the solution to the following problem:\\
\noindent
\begin{align*}
\text{min}    &\quad  \ell(c) \\
\text{subject to} & \quad \birth(c)=b \\
                  & \quad \death(c)=d\\
                   & \quad c \in \ker\left (\bd_{p}\big|_{C_p(K_b)}\right )
\end{align*}
\noindent
This problem was first studied in \cite{CF08}.

\section{Time-optimal representatives}

We start by discussing an example that will guide us in finding a suitable notion of time-optimal PH cycle representative. Consider the  simplicial complex associated to the point cloud illustrated in Figure \ref{F:mot ex}(a), which we can think of as being obtained by a time-delay embedding of a noisy sine curve.

Among the different possible choices of  $1$-cycle representatives, we give four examples in Figure \ref{F:mot ex}(b), highlighted in green.
Which of these choices should we consider as optimal with respect to time?
For the  problem we are interested in studying in this paper, namely, studying univariate time series through time-delay embeddings, 
we are interested in obtaining cycle representatives that correspond to time-series values that are not too far from each other in time, since this allows us to interpret  the topological feature (i.e., a non-trivial PH class), through a selection of data points in the original time series that correspond to that feature. For instance, for the example application that we study in this paper --- delayed oscillator models of the El Nin\~o Southern Oscillation ---, this allows one to ask the question of whether the topological features that one recovers have any physical meaning at all.

In particular, for the noisy sine curve from Figure \ref{F:mot ex}(a), this means that we would want to choose representative $1$-cycles whose vertices correspond to time-series values that are contained in an interval of length $2\pi$, the period of the sine curve. Thus, 
we make the following observations:
\begin{itemize}
\item Any reasonable notion of time-optimal cycle representative should not choose a $1$-cycle that includes the cyan edges from Figure~\ref{F:mot ex}(a).  
\item On the other hand, it is likely that we will have to make a choice between the edges incident to vertices labelled by $0$ and $\pi/3$, $5/3\pi$ and those incident to vertices labelled by $2\pi$ and $\pi/3$, $5/3\pi$.
\end{itemize}
Thus, based on these observations, we would wish to consider the two cycles on the left of Figure~\ref{F:mot ex, repr} as time-optimal, but not the two ones on the right.

We  can thus formulate our problem as follows. We let $f\colon \mathbb{R}\times \mathbb{R}\to \mathbb{R}_{\geq 0 }$ be a norm. To each edge $e=\{s,s'\}$ we assign the weight $w(e):=f(T(s),T(s'))$, where for a point $s$ in the point cloud, we denote by $T(s)$ its time label. Then for each possible $1$-cycle representative, we compute its cost by considering how far apart its edges are in time. More precisely, to compute the cost of  each edge, we compute the difference between its weight and the weights of the edges adjacent to it, and we retain the maximum of these differences as the cost of the edge. The cost of the $1$-cycle is then the sum over the costs of its edges. 
A time-optimal cycle is  one that minimises this cost function.

We can thus write our optimisation problem as follows:

$$
c_{\mathrm{opt}}=\underset{c \text{ a } 1-cycle}{\argmin}\left \{ \sum_{e\in c}\underset{\text{adjacent to }e}{\underset{e'\in c}{\max}} \big|w(e)-w(e')\big| \right\} \, .
$$

\begin{figure}[h!]

\begin{subfigure}{0.4\textwidth}
\begin{tikzpicture}
\node[label={$4/3\pi$}] (4/3pi) at (0,0) {$\bullet$};
\node[label={$\pi$}] (pi) at (2,0.6) {$\bullet$};
\node[label=below:{$3\pi$}] (3pi) at (2.45,0.1) {$\bullet$};
\node[label={$2/3\pi$}] (2/3pi) at (2.8,3) {$\bullet$};
\node[label={$\pi/3$}] (pi/3) at (0.8,5) {$\bullet$};
\node[label={[shift={(0.3,-0.7)}]$0$}] (0) at (-1.8,3.8) {$\bullet$};
\node[label={$2\pi$}] (2pi) at (-2.4,4.2) {$\bullet$};
\node[label=below:{$5/3\pi$}] (5/3pi) at (-2,1.2) {$\bullet$};
\draw[-] (pi)--(3pi);
 \draw[-] (4/3pi)--(pi)--(2/3pi)--(pi/3)--(0)--(5/3pi)--(4/3pi);
 \draw[-,cyan] (4/3pi)--(3pi)--(2/3pi);
 \draw[-,orange] (5/3pi)--(2pi)--(pi/3);
 \draw[-](0)--(2pi);
  \begin{pgfonlayer}{bg}    
 \fill [gray!10] (4/3pi.center) -- (pi.center) -- (3pi.center) -- cycle;
 \fill [gray!10] (2/3pi.center) -- (pi.center) -- (3pi.center) -- cycle;
 \fill [gray!10] (2pi.center) -- (0.center) -- (5/3pi.center) -- cycle;
 \fill [gray!10] (2pi.center) -- (0.center) -- (pi/3.center) -- cycle;
   \end{pgfonlayer}
\end{tikzpicture}
\subcaption{A simplicial complex with time-labelled vertices, associated to a point cloud arising from a noisy sine curve.}
\label{F:mot ex, simp com}
\end{subfigure}\qquad \qquad
\begin{subfigure}{0.5\textwidth}
\begin{tikzpicture}[scale=0.48, transform shape]
\node[label={$4/3\pi$}] (4/3pi) at (0,0) {$\bullet$};
\node[label={$\pi$}] (pi) at (2,0.6) {$\bullet$};
\node[label=below:{$3\pi$}] (3pi) at (2.45,0.1) {$\bullet$};
\node[label={$2/3\pi$}] (2/3pi) at (2.8,3) {$\bullet$};
\node[label={$\pi/3$}] (pi/3) at (0.8,5) {$\bullet$};
\node[label={[shift={(0.3,-0.7)}]$0$}] (0) at (-1.8,3.8) {$\bullet$};
\node[label={$2\pi$}] (2pi) at (-2.4,4.2) {$\bullet$};
\node[label=below:{$5/3\pi$}] (5/3pi) at (-2,1.2) {$\bullet$};
\draw[-] (pi)--(3pi);
 \draw[-] (4/3pi)--(pi)--(2/3pi)--(pi/3)--(0)--(5/3pi)--(4/3pi);
 \draw[-,black] (4/3pi)--(3pi)--(2/3pi);
 \draw[-,black] (5/3pi)--(2pi)--(pi/3);
 \draw[-](0)--(2pi);
  \draw[-,thick, green] (4/3pi)--(pi)--(2/3pi)--(pi/3)--(2pi)--(5/3pi)--(4/3pi);
  \begin{pgfonlayer}{bg}    
 \fill [gray!10] (4/3pi.center) -- (pi.center) -- (3pi.center) -- cycle;
 \fill [gray!10] (2/3pi.center) -- (pi.center) -- (3pi.center) -- cycle;
 \fill [gray!10] (2pi.center) -- (0.center) -- (5/3pi.center) -- cycle;
 \fill [gray!10] (2pi.center) -- (0.center) -- (pi/3.center) -- cycle;
   \end{pgfonlayer}
\end{tikzpicture}\qquad
\begin{tikzpicture}[scale=0.48, transform shape]
\node[label={$4/3\pi$}] (4/3pi) at (0,0) {$\bullet$};
\node[label={$\pi$}] (pi) at (2,0.6) {$\bullet$};
\node[label=below:{$3\pi$}] (3pi) at (2.45,0.1) {$\bullet$};
\node[label={$2/3\pi$}] (2/3pi) at (2.8,3) {$\bullet$};
\node[label={$\pi/3$}] (pi/3) at (0.8,5) {$\bullet$};
\node[label={[shift={(0.3,-0.7)}]$0$}] (0) at (-1.8,3.8) {$\bullet$};
\node[label={$2\pi$}] (2pi) at (-2.4,4.2) {$\bullet$};
\node[label=below:{$5/3\pi$}] (5/3pi) at (-2,1.2) {$\bullet$};
\draw[-] (pi)--(3pi);
 \draw[-] (4/3pi)--(pi)--(2/3pi)--(pi/3)--(0)--(5/3pi)--(4/3pi);
 \draw[-,black] (4/3pi)--(3pi)--(2/3pi);
 \draw[-,black] (5/3pi)--(2pi)--(pi/3);
 \draw[-](0)--(2pi);
 \draw[-,thick, green] (4/3pi)--(3pi)--(2/3pi)--(pi/3)--(0)--(5/3pi)--(4/3pi);
  \begin{pgfonlayer}{bg}    
 \fill [gray!10] (4/3pi.center) -- (pi.center) -- (3pi.center) -- cycle;
 \fill [gray!10] (2/3pi.center) -- (pi.center) -- (3pi.center) -- cycle;
 \fill [gray!10] (2pi.center) -- (0.center) -- (5/3pi.center) -- cycle;
 \fill [gray!10] (2pi.center) -- (0.center) -- (pi/3.center) -- cycle;
   \end{pgfonlayer}
\end{tikzpicture}
\begin{tikzpicture}[scale=0.48, transform shape]
\node[label={$4/3\pi$}] (4/3pi) at (0,0) {$\bullet$};
\node[label={$\pi$}] (pi) at (2,0.6) {$\bullet$};
\node[label=below:{$3\pi$}] (3pi) at (2.45,0.1) {$\bullet$};
\node[label={$2/3\pi$}] (2/3pi) at (2.8,3) {$\bullet$};
\node[label={$\pi/3$}] (pi/3) at (0.8,5) {$\bullet$};
\node[label={[shift={(0.3,-0.7)}]$0$}] (0) at (-1.8,3.8) {$\bullet$};
\node[label={$2\pi$}] (2pi) at (-2.4,4.2) {$\bullet$};
\node[label=below:{$5/3\pi$}] (5/3pi) at (-2,1.2) {$\bullet$};
\draw[-] (pi)--(3pi);
 \draw[-] (4/3pi)--(pi)--(2/3pi)--(pi/3)--(0)--(5/3pi)--(4/3pi);
 \draw[-,black] (4/3pi)--(3pi)--(2/3pi);
 \draw[-,black] (5/3pi)--(2pi)--(pi/3);
 \draw[-](0)--(2pi);
 \draw[-,thick, green] (4/3pi)--(pi)--(2/3pi)--(pi/3)--(0)--(5/3pi)--(4/3pi);
  \begin{pgfonlayer}{bg}    
 \fill [gray!10] (4/3pi.center) -- (pi.center) -- (3pi.center) -- cycle;
 \fill [gray!10] (2/3pi.center) -- (pi.center) -- (3pi.center) -- cycle;
 \fill [gray!10] (2pi.center) -- (0.center) -- (5/3pi.center) -- cycle;
 \fill [gray!10] (2pi.center) -- (0.center) -- (pi/3.center) -- cycle;
   \end{pgfonlayer}
\end{tikzpicture}\qquad
\begin{tikzpicture}[scale=0.48, transform shape]
\node[label={$4/3\pi$}] (4/3pi) at (0,0) {$\bullet$};
\node[label={$\pi$}] (pi) at (2,0.6) {$\bullet$};
\node[label=below:{$3\pi$}] (3pi) at (2.45,0.1) {$\bullet$};
\node[label={$2/3\pi$}] (2/3pi) at (2.8,3) {$\bullet$};
\node[label={$\pi/3$}] (pi/3) at (0.8,5) {$\bullet$};
\node[label={[shift={(0.3,-0.7)}]$0$}] (0) at (-1.8,3.8) {$\bullet$};
\node[label={$2\pi$}] (2pi) at (-2.4,4.2) {$\bullet$};
\node[label=below:{$5/3\pi$}] (5/3pi) at (-2,1.2) {$\bullet$};
\draw[-] (pi)--(3pi);
 \draw[-] (4/3pi)--(pi)--(2/3pi)--(pi/3)--(0)--(5/3pi)--(4/3pi);
 \draw[-,black] (4/3pi)--(3pi)--(2/3pi);
 \draw[-,black] (5/3pi)--(2pi)--(pi/3);
 \draw[-](0)--(2pi);
 \draw[-,thick, green] (4/3pi)--(pi)--(2/3pi)--(pi/3)--(0)--(2pi)--(5/3pi)--(4/3pi);
  \begin{pgfonlayer}{bg}    
 \fill [gray!10] (4/3pi.center) -- (pi.center) -- (3pi.center) -- cycle;
 \fill [gray!10] (2/3pi.center) -- (pi.center) -- (3pi.center) -- cycle;
 \fill [gray!10] (2pi.center) -- (0.center) -- (5/3pi.center) -- cycle;
 \fill [gray!10] (2pi.center) -- (0.center) -- (pi/3.center) -- cycle;
   \end{pgfonlayer}
\end{tikzpicture}
\subcaption{Different representative $1$-cycles (in green).}
\label{F:mot ex, repr}
\end{subfigure}

\caption{}\label{F:mot ex}
\end{figure}

While in our example we have considered $1$-cycles, one can in a similar way define an optimization problem for higher-dimensional cycles. In what follows our discussion will not be restricted to any particular  cycle dimension.

\subsection{Loss function}
We now look at defining an appropriate loss function for the minimization problem.
Let $c \in C_p(K)$ be a chain. We wish to build a non-negative matrix $W$ that appropriately weighs the chain $c$ such that the optimal solution is a chain with minimal dispersion in time. 
We first  make the notion of ``time-closeness'' or ``minimal time dispersion'' precise:

\begin{definition}
Given a $p$-chain, its \define{time dispersion} is the difference between the minimum and maximum time labels of any of the vertices contained in the $p$-simplices in the chain.
\end{definition}

We note that while for other types of applications, e.g., temporal networks (see also discussion in Section \ref{S: conclusion}), one might be interested in considering another notions of time dispersion, the notion of time dispersion we consider here is motivated by the study of time series, for which we want to obtain data points in the time series that are not too far apart, or dispersed.

We next choose an order on the simplices in the chain $c$, and label the simplices $\sigma_1,\dots , \sigma_n$ according to this order. 
Each entry $w_{ii}$ of the matrix describes the 
cost
of selecting the $p$-simplex $\sigma_i$ in the chain.
Similarly, the entry $w_{ij}$ represents the cost
of having both $\sigma_i$ and $\sigma_j$ as part of the chain.
Our goal is then to minimize the norm  $\lVert Wc \rVert$.
In the following we consider two different ways of measuring time-optimality:
\begin{enumerate}
    \item ``vertex-based'': We find chains with vertices with  time labels that are close to each other.
    \item ``simplex-based'': We first situate each simplex in time, and find a chain of closely situated simplices.
\end{enumerate}

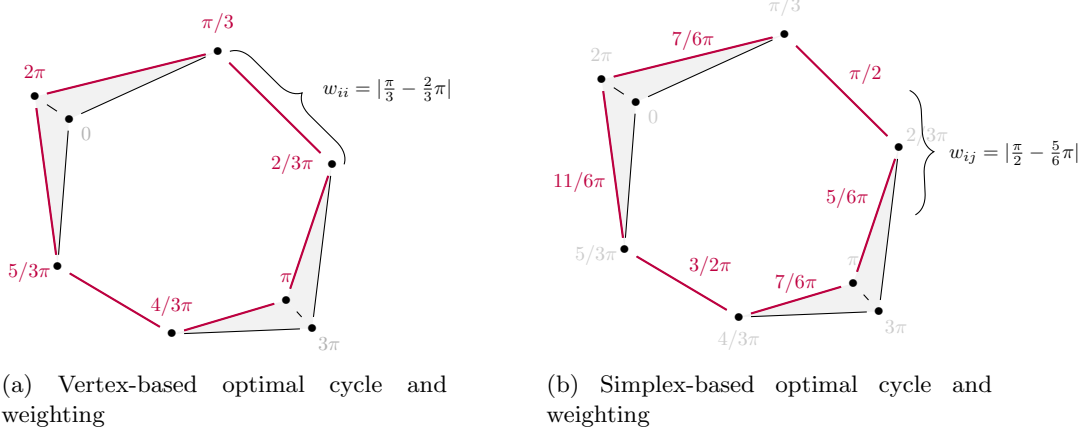
\begin{figure}
    \centering
    \begin{subfigure}{0.45\textwidth}
    \centering
    \begin{tikzpicture}[scale=0.75, transform shape]
    \node[label={[scale=1,purple]$4/3\pi$}] (4/3pi) at (0,0) {$\bullet$};
    \node[label={[scale=1,purple]$\pi$}] (pi) at (2,0.6) {$\bullet$};
    \node[label={[scale=1, black!30,shift={(0.3,-0.7)}]$3\pi$}] (3pi) at (2.45,0.1) {$\bullet$};
    \node[label={[scale=1,purple, left, shift={(-0.2,-0.2)}]$2/3\pi$}] (2/3pi) at (2.8,3) {$\bullet$};
    \node[label={[scale=1,purple]$\pi/3$}] (pi/3) at (0.8,5) {$\bullet$};
    \node[label={[scale=1,black!30,shift={(0.3,-0.7)}]$0$}] (0) at (-1.8,3.8) {$\bullet$};
    \node[label={[scale=1,purple]$2\pi$}] (2pi) at (-2.4,4.2) {$\bullet$};
    \node[label={[below left,scale=1,purple]$5/3\pi$}] (5/3pi) at (-2,1.2) {$\bullet$};
    \draw[-] (pi)--(3pi);
    \draw[-] (4/3pi)--(pi)--(2/3pi)--(pi/3)--(0)--(5/3pi)--(4/3pi);
    \draw[-,black] (4/3pi)--(3pi)--(2/3pi);
    \draw[-,black] (5/3pi)--(2pi)--(pi/3);
    \draw (0)--(2pi);
    \draw[-, thick, purple] (4/3pi) -- (pi);
    \draw[-, thick, purple] (pi) -- (2/3pi);
    \draw[-, thick, purple] (2/3pi) -- (pi/3);
    \draw[-, thick, purple] (pi/3) -- (2pi);
    \draw[-, thick, purple] (2pi) -- (5/3pi);
    \draw[-, thick, purple] (5/3pi) -- (4/3pi);
    
    \draw[decorate, decoration={brace, amplitude=8pt}] 
        (pi/3.east) -- (2/3pi.east) node[midway, above, xshift=50pt, scale=1] 
        {$w_{ii} = |\frac{\pi}{3} - \frac{2}{3}\pi|$};
    \begin{pgfonlayer}{bg}
    \fill [gray!10] (4/3pi.center) -- (pi.center) -- (3pi.center) -- cycle;
    \fill [gray!10] (2/3pi.center) -- (pi.center) -- (3pi.center) -- cycle;
    \fill [gray!10] (2pi.center) -- (0.center) -- (5/3pi.center) -- cycle;
    \fill [gray!10] (2pi.center) -- (0.center) -- (pi/3.center) -- cycle;
    \end{pgfonlayer}
    \end{tikzpicture}
     \subcaption{Vertex-based optimal cycle and weighting}
    \label{fig:vertex-based-illustration}
    \end{subfigure}
    \hfill
    \begin{subfigure}{0.45\textwidth}
        \centering
        \begin{tikzpicture}[scale=0.75, transform shape]
        \node[label={[scale=1,below,black!20,shift={(0.0,-0.3)}]$4/3\pi$}] (4/3pi) at (0,0) {$\bullet$};
        \node[label={[scale=1,black!20]$\pi$}] (pi) at (2,0.6) {$\bullet$};
        \node[label={[scale=1,black!20, shift={(0.3,-0.7)}]$3\pi$}] (3pi) at (2.45,0.1) {$\bullet$};
        \node[label={[scale=1,black!20,right]$2/3\pi$}] (2/3pi) at (2.8,3) {$\bullet$};
        \node[label={[scale=1,black!20]$\pi/3$}] (pi/3) at (0.8,5) {$\bullet$};
        \node[label={[scale=1,black!20,shift={(0.3,-0.7)}]$0$}] (0) at (-1.8,3.8) {$\bullet$};
        \node[label={[scale=1,black!20]$2\pi$}] (2pi) at (-2.4,4.2) {$\bullet$};
        \node[label={[below left,scale=1,black!20]$5/3\pi$}] (5/3pi) at (-2,1.2) {$\bullet$};
        \draw[-] (pi)--(3pi);
        \draw[-] (4/3pi)--(pi)--(2/3pi)--(pi/3)--(0)--(5/3pi)--(4/3pi);
        \draw[-,black] (4/3pi)--(3pi)--(2/3pi);
        \draw[-,black] (5/3pi)--(2pi)--(pi/3);
        \draw[-] (0)--(2pi);
        \draw[-, thick, purple] (4/3pi) -- node[midway, above, scale=1]{$7/6\pi$} (pi);
        \draw[-, thick, purple] (pi) -- node[midway, above left, scale=1]{$5/6\pi$} (2/3pi);
        \draw[-, thick, purple] (2/3pi) -- node[midway, above right, scale=1]{$\pi/2$} (pi/3);
        \draw[-, thick, purple] (pi/3) -- node[midway, above, scale=1]{$7/6\pi$} (2pi);
        \draw[-, thick, purple] (2pi) -- node[midway, below left, scale=1]{$11/6\pi$} (5/3pi);
        \draw[-, thick, purple] (5/3pi) -- node[midway, above right, scale=1]{$3/2\pi$} (4/3pi);

        \draw[decorate,
          decoration={brace,
                       amplitude=8pt,
                       raise=4pt,
                       aspect=0.5,
                       mirror}] 
         let \p1 = ($(pi)!0.5!(2/3pi) + (0.5, 0)$),
             \p2 = ($(2/3pi)!0.5!(pi/3)  + (1, 0)$)
         in (\p1) -- node[midway,right,scale=1,xshift=20pt] {$w_{ij} = |\frac{\pi}{2}-\frac{5}{6}\pi|$} (\p2);

        \begin{pgfonlayer}{bg}
        \fill [gray!10] (4/3pi.center) -- (pi.center) -- (3pi.center) -- cycle;
        \fill [gray!10] (2/3pi.center) -- (pi.center) -- (3pi.center) -- cycle;
        \fill [gray!10] (2pi.center) -- (0.center) -- (5/3pi.center) -- cycle;
        \fill [gray!10] (2pi.center) -- (0.center) -- (pi/3.center) -- cycle;
        \end{pgfonlayer}
        \end{tikzpicture}
        \subcaption{Simplex-based optimal cycle and weighting}
        \label{fig:simplex-based-illustration}
    \end{subfigure}%
    \caption{Illustrations of the difference between considering simplex weighting and vertex weighting for the example from Figure \ref{F:mot ex, repr}. In both situations the highlighted cycle represents the optimal solution. \textbf{a)} In vertex-based weighting,  the cost $w_{ii}$ of selecting  1-simplex $\sigma_i$ is given by the difference of the time labels of its vertices. \textbf{b)} For simplex-based weighting, we first assign a time label to each simplex,  as the mean of the time label of its vertices. The cost of a chain is then the sum of the differences $w_{ij}$ of the time labels of adjacent simplices $\sigma_i,\sigma_j$.  }
    \label{F:mot exp time-opt repr}
\end{figure}

\subsubsection{Vertex-based time optimality}
In this scenario the cost of selecting a $p$-simplex depends only on the time labels of its vertices. 
Consider a $1$-simplex $\sigma = \{v_i,v_j\}$. If we are attempting to find a $1$-cycle that minimizes the dispersion in time the intuitive cost of selecting simplex $\sigma$ would be the difference in the time label of the vertices $|T(v_i) - T(v_j)|$.
The optimal solution is a $1$-cycle composed of $1$-simplices with vertices that are as close  in time as possible.

Given a $p$-simplex $\sigma_i = \{v_{i_0}, v_{i_1},\dots,v_{i_{p+1}}\}$  we let 
\[
{T_{i}^{\max} = \max \left \{T(v_{i_0}), T(v_{i_1}),\dots,T(v_{i_{p+1}})\right \}}
\] be the maximum time label of the vertices of $\sigma_i$.
Similarly, we let $T_{i}^{\min}$ be the minimum time label of the vertices of $\sigma_i$.
Then we define the weight matrix $W$ as the diagonal matrix with diagonal entries given by:
$$
w_{ii} = 
    T_{i}^{\max} - T_{i}^{\min}\, .
$$
We thus have that the cost of selecting a chain containing $\sigma_i$ is given by the maximum difference of the time labels of its vertices.

\subsubsection{Simplex-based time optimality}\label{S:simplex-based time opt}
In a simplex-based approach we first assign a time label $T(\sigma)$ to each $p-$simplex.
This can be seen as positioning the simplex $\sigma$ in time.
The solution of the optimization problem is a cycle composed of simplices where adjacent simplices are closely placed in time.
\begin{definition}
Two $p$-simplices $\sigma_i, \sigma_j$ are called \define{adjacent} if their intersection is non-empty and is of cardinality $p$. That is, if there exists a $p-1$-simplex $\tau$ such that $\sigma_i \cap \sigma_j = \tau$.
\end{definition}

Thus, we consider the weight matrix with entries defined as follows:
$$
w_{ij} = \begin{cases}
    |T(\sigma_i)- T(\sigma_j)| &\quad \text{if } \sigma_i, \sigma_j \text{ are adjacent}\\
    0 &\quad \text{otherwise}
\end{cases}
$$
In our experiments we consider the time label of a $p$-simplex $\sigma =\{v_{i_0}, v_{i_1},\dots,v_{i_{p+1}}\}$ as the mean time label of its vertices: 
$$
T(\sigma) = \frac{1}{(p+1)}\sum^{p+1}_{k=0}T(v_{i_k})\, .
$$

In Figure \ref{F:mot exp time-opt repr} we illustrate the two notions of vertex- and simplex-optimal $1$-cycle representatives obtained for the example in Figure \ref{F:mot ex, repr}.
For practical applications, we need to further modify our optimization problem, as we explain next in Section \ref{A:approx}.

\subsection{Approximate PH cycle representatives}\label{A:approx}
When searching for a representative cycle of a PH  class one forces the solution to have the same birth and death values $(b,d)$ as the class it represents.
In other words, a solution must contain the birth simplex and cannot contain any simplex that appears after the birth value.
As we find in addressing our problem, this constraint limits the possible solution set of the optimization problem in a way that is too restrictive to give any meaningful cycles.

We thus  relax this constraint and allow the representative cycle to not necessarily have the same persistence  of the class it represents. 
More precisely, given a minimum persistence value $\varepsilon$ we instead search for representatives with a birth time of $(d -\varepsilon )$, that is, with a persistence of $\varepsilon$.

\begin{figure}[b!]
 \centering
\includegraphics[width=0.8\textwidth]{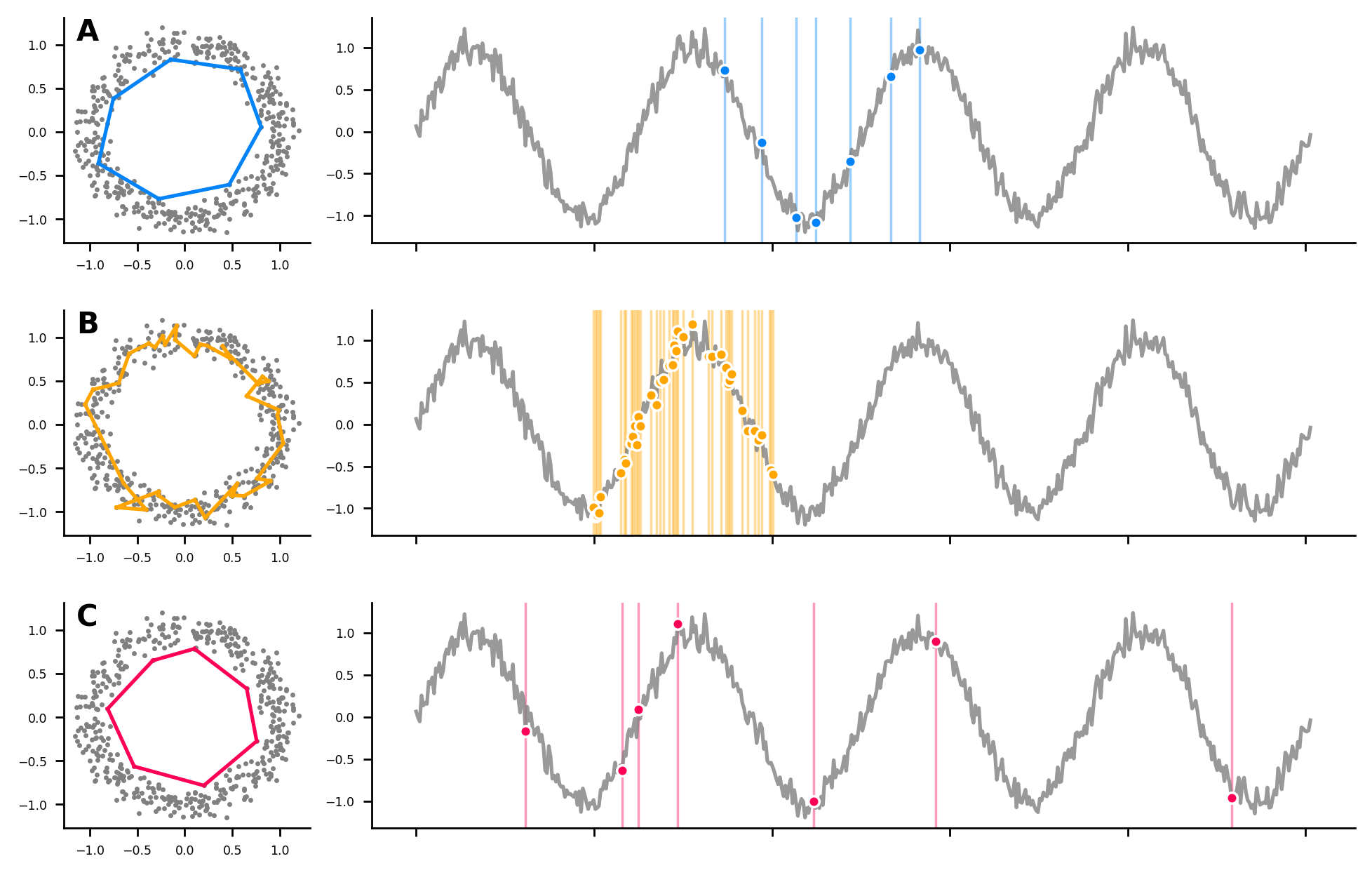}
\caption{We illustrate different $1$-cycle PH representatives for the optimal embedding of $f(t) = \sin{t}+\epsilon$, where $\epsilon$ is a random variable with normal distribution,  in the point cloud (left), as well as in the original time series (right). \textbf{A}: Simplex-based time optimal cycle. \textbf{B}: Vertex-based time optimal cycle. \textbf{C}: Length-optimal cycle. We note that both simplex and vertex-based optimal cycles have time dispersion that is close to the period of the sine, and thus close to the desired minimum. On the other hand, the length-optimal cycle is much further spread out in time, and thus doesn't provide a meaningful cycle representative.}\label{F:optimal cycles sine}

\end{figure}

One could then ask what a good choice for such a minimum  persistence value $\varepsilon$ might be. In our set up we follow the computational methodology from \cite{GP24}, which provides lower bounds on persistence that distinguish between topological noise (points in the persistence diagram with persistence smaller than the bound) and topological signal (points in the persistence diagram with persistence greater or equal than the bound). Thus, we take use this bound for several of the examples considered 
(in the noisy sine example, see Figure \ref{F:optimal cycles sine} and the double sine, see Figure \ref{fig:example-doublesine-embedding}). For computational reasons, for the time series arising from the ENSO models we take instead $90\%$ of the persistence of the PH class (which thus gives a bound smaller than the significance bound from \cite{GP24}), namely, we take
 $\varepsilon = 0.9(d-b)$ (see Figures \ref{F:ENSO 1.65},\ref{F:ENSO 1.4} and \ref{F:ENSO 1.9}). 

\begin{figure}[t!]
\centering
\includegraphics[width=0.8\textwidth]{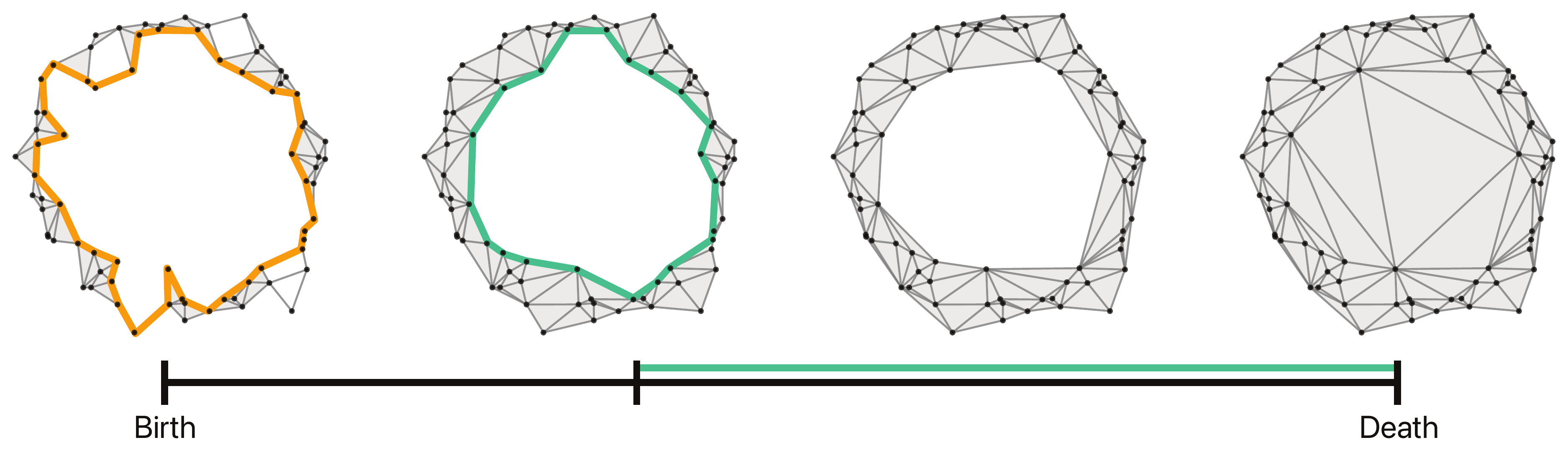}
\caption{Example of representative $1$-cycles with minimal $l_1$ norm for different persistence values (see Appendix \ref{A:lin progr} for details). Left/orange: we depict a $1$-cycle with full persistence. Second from left/green: we depict a $1$-cycle homologous to the previous one, but with smaller persistence, and with smaller $l_1$ norm.  We thus note that relaxing the persistence of the representative allows for a solution with smaller $l_1$ norm.}
\label{fig:relaxation illustration}
\end{figure}

We note that, as we have already observed in Section \ref{SS:PH}, given filtration values $r\leq r'$, we have that $K_r\subset K_{r'}$ and in particular, we can identify $C_p(K_r)$ with a vector subspace of $C_p(K_{r'})$.
Therefore, we have that any chain that exists in $C_p(K_r)$ also exists in $C_p(K_{r'})$.
More specifically, given any optimization function $f$ we have that:
$$
r\leq r' \implies \min \{f(c) \mid c \in C_p(K_r)\} \geq \min \{f(c) \mid c \in C_p(K_{r'})\}\, .
$$

In practice, if we allow the representative cycle to have a smaller persistence value than that of its class, we are guaranteed to obtain a cycle with a smaller loss value. 
This can very beneficial, particularly in the following situations: 
\begin{enumerate}
    \item One has a lower bound on  persistence, allowing to distinguish between topological ``signal'' and ``noise''. This can be given for instance, as in our case, by the computational methodology from \cite{GP24} to compute embeddings of quasi-periodic time series. More generally, such a bound may be obtained as a result of a suitable statistical analysis. In all such cases, one could then use such a lower bound as a lower bound for the persistence of the cycle representative. 
    \item If it's likely that a birth simplex may be suboptimal for the application at hand. For instance, in an application such as the one discussed in this paper, for which one seeks time-optimal cycles, it might happen that the birth simplex connects vertices very far apart in time, and is thus  a  poor choice. 
    In such a case, even small relaxations, for instance, taking a representative with $95\%$ of the persistence of its class, can  lead to much better outcomes, since the solution to the optimization problem  no longer has to include the given birth simplex.
\end{enumerate}

Given a homology class with barcode $(b,d)$ and an initial representative cycle $c_0$, searching for an optimal homologous cycle with minimum persistence $\varepsilon$ involves solving the problem:
\begin{align*}
\text{min}   & \quad \ell(c) \\
\text{subject to} & \quad c = c_0 + \bd_{p+1}(w) \\
                  & \quad w \in C_{p+1}(K_{d-\varepsilon})
\end{align*}

In Figure \ref{F:optimal cycles sine} we illustrate how our notions of time-optimal PH cycles compare with  the existing notion of length-optimal PH cycles on a synthetic example of a noisy sine curve. 
We provide details about the algorithms and implementation in Appendix \ref{A:details algo impl}.

\section{Experiments}
In all experiments we consider a Vietoris-Rips filtration on the embedded time-series (see Appendix \ref{A:optimal par}).
For each PH class we obtain an initial representative $\mathbf{c}_0$ using the $R=\bd \cdot V$ decomposition of the filtration boundary matrix $\bd$.
For each PH class with birth and death values $(b,d)$ we restrict the domain and codomain of the boundary operator by considering the sets:
$$
P = \{\sigma \in S_{p}(K_b) \mid \text{birth}(\sigma) \leq b\}\\
$$
$$
\hat{Q} = \{\sigma \in S_{p+1}(K_b) \mid \text{birth}(\sigma) \leq b \quad \text{and}\quad R[:,\sigma] \neq 0 \}\, ,\\
$$
which corresponds to the set of $p$- and $p+1$-simplices that are alive at filtration time $b$. This assures that the solution has a persistence value of at least $(d-b)$ (see Appendix \ref{A:details algo impl} for more details). For a relaxed problem with minimum persistence $\varepsilon$ we simply take $b'=d-\varepsilon$.
We then solve the following linear problem \cite{escolar2016optimal,lu2021minimal}:
\begin{align*}
\text{min}   &  \quad \lVert  W \mathbf{c} \rVert_1 = \sum_i \sum_j w_{ij} (c_{j}^+ + c_{j}^-) \\
\text{subject to  } & (\mathbf{c}^+ - \mathbf{c}^-) = \mathbf{c_0} + \bd_{p+1}[P,\hat{Q}](\mathbf{w}) \\
            & \mathbf{w}\in \mathbb{R}^{|\hat{Q}|}\\
            & \mathbf{c} \in \mathbb{R}^{|P|}\\
            & \mathbf{c}^+,\mathbf{c}^- \geq 0
\end{align*}
where $\bd_{p+1}[P,\hat{Q}]$ is a submatrix obtained by selecting the rows and columns of the boundary matrix $\bd_{p+1}$.

The experiments were run using the Gurobi solver \cite{gurobi}  on an Apple M1 Pro with 16Gb of memory.
The code to reproduce our experiments is available at our \href{https://github.com/antonio-leitao/optimal-cycles}{GitHub repository} \cite{github}.

\subsection{Synthetic quasi-periodic time series}

We consider the quasi-periodic time series given by the function

\begin{equation}\notag
    f(t) = 2\sin{t} + 1.8sin{\sqrt{3}t}\quad 0\leq t \leq 60\pi  \, .
\end{equation}

\begin{figure}[t!]
\centering
\begin{subfigure}[b]{0.45\textwidth}
  \centering
\includegraphics[height=5cm,keepaspectratio]{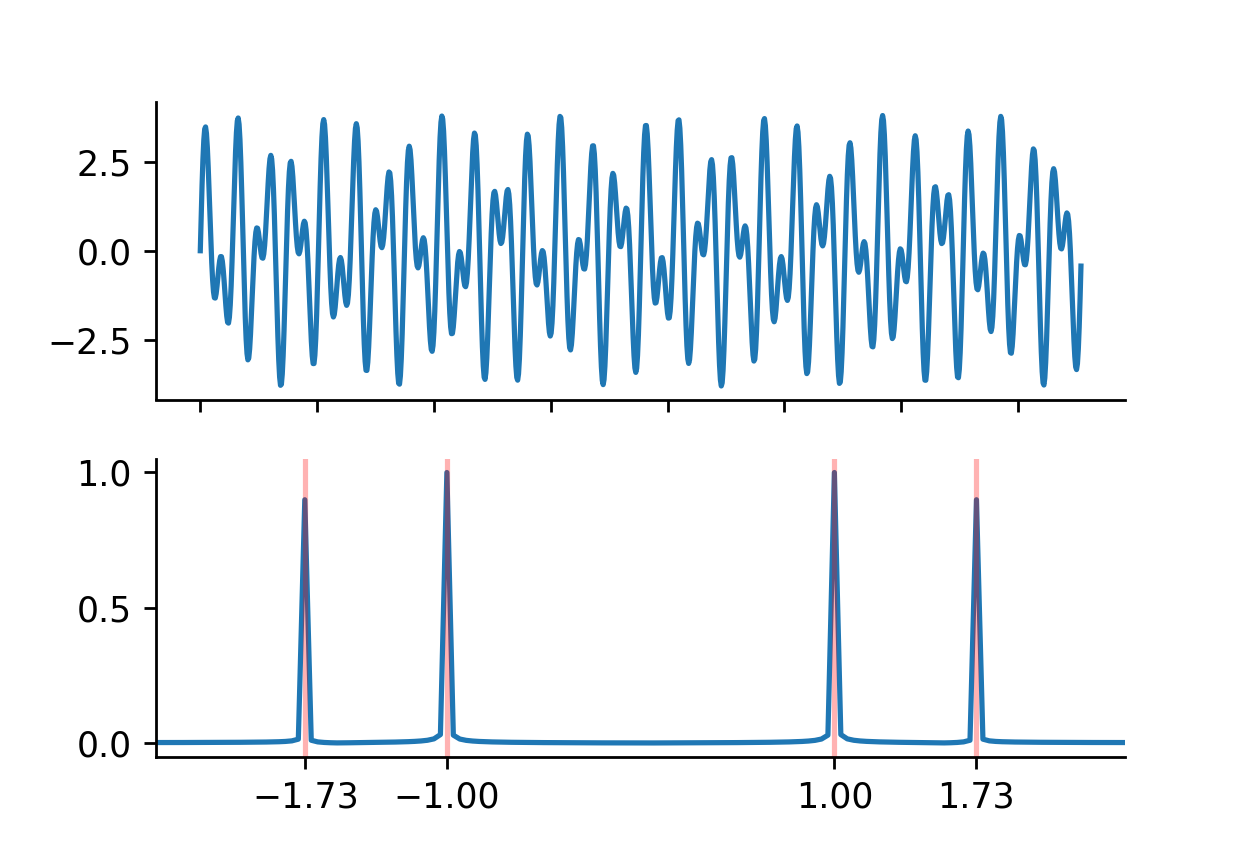}
  \subcaption{Original Signal}
\end{subfigure}%
\begin{subfigure}[b]{0.4\textwidth}
  \centering
 \includegraphics[height=5cm,keepaspectratio]{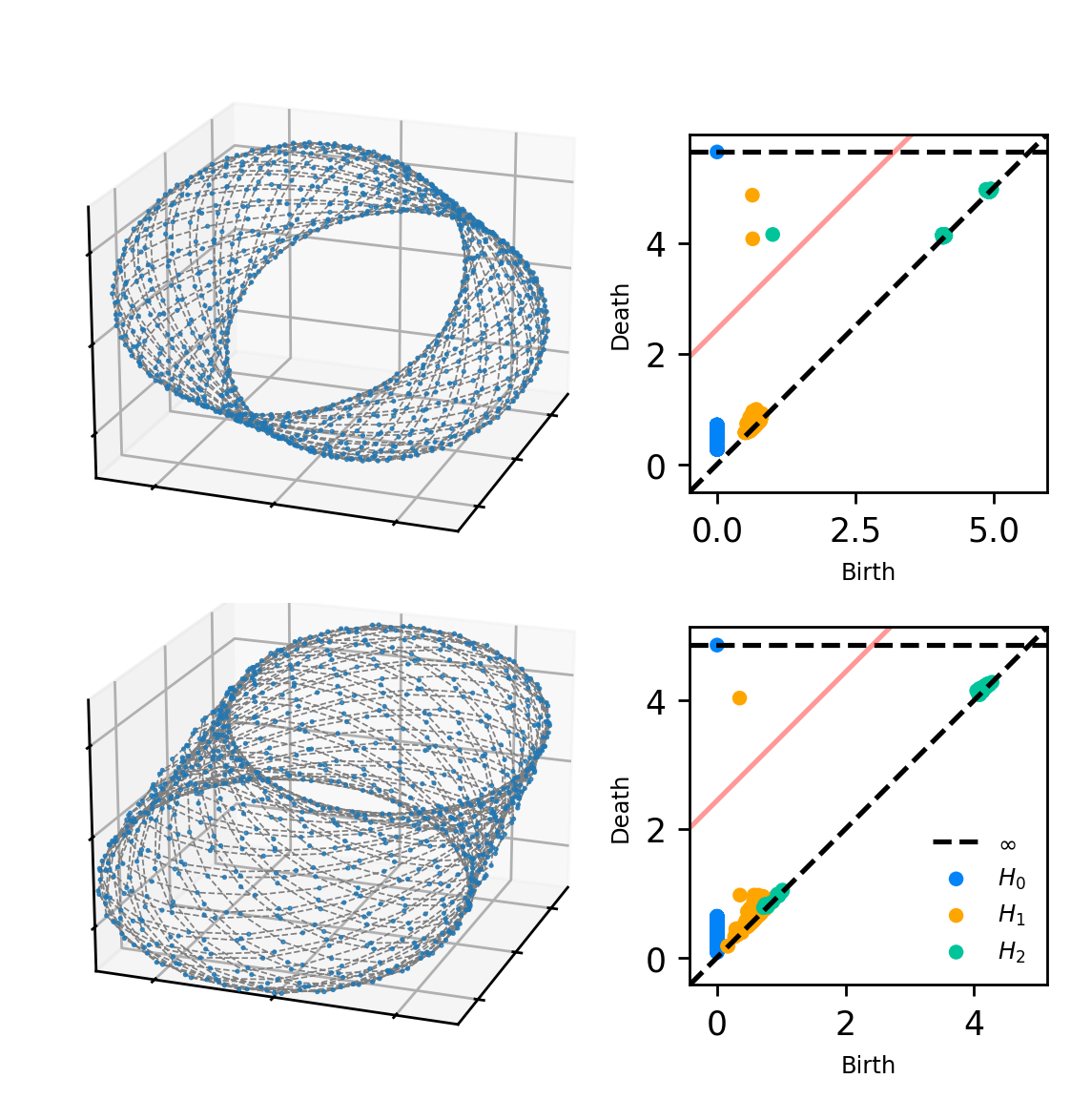} 
  \subcaption{Sliding window embedding}
\end{subfigure}

\caption{Pipeline for computing the optimal embedding parameters for $f(t) = 2\sin{t} + 1.8sin{\sqrt{3}}t$. \textbf{a)} Top: graph of the function. Bottom: Fourier diagram, with peaks at $1$ and $\sqrt{3}$. \textbf{b)} Two different embeddings and corresponding PDs for two different values of the delay $\tau$.
Top left: embedding for optimal value of the delay $\tau$.
Top right: corresponding PD.
The continuous red diagonal line indicates the lower bound on persistence for all homology degrees (see Appendix \ref{A:persistence bounds} for details).
Bottom: embedding for suboptimal delay value, together with the corresponding PD. The embeddings are in  $\mathbb{R}^4$ and we show the projection onto the first three principal components.}
\label{fig:example-doublesine-embedding}
\end{figure}

We sample $1000$ points and compute the optimal embedding parameters using the methodology introduced in \cite{perea2019topological}, see Figure \ref{fig:example-doublesine-embedding} for an illustration of the pipeline, and Appendix \ref{A:optimal par} for details. 
We then proceed to compute vertex- and simplex-based time-optimal PH $1$-cycle representatives for the embedding corresponding to the optimal parameters (top embedding in Figure \ref{fig:example-doublesine-embedding}(b)).

Given a  a lower bound on the persistence of  PH classes  that we may interpret as ``significant'' \ref{fig:example-doublesine-embedding}(b), we use that same lower bound as a minimum persistence value for the  search of approximate representatives (as opposed to representatives with same persistence as the  class). The representative optimization is done with half the points (every other point sampled) due to computational constraints.

 \begin{figure}[h!]
 \centering
\includegraphics[width=0.8\textwidth]{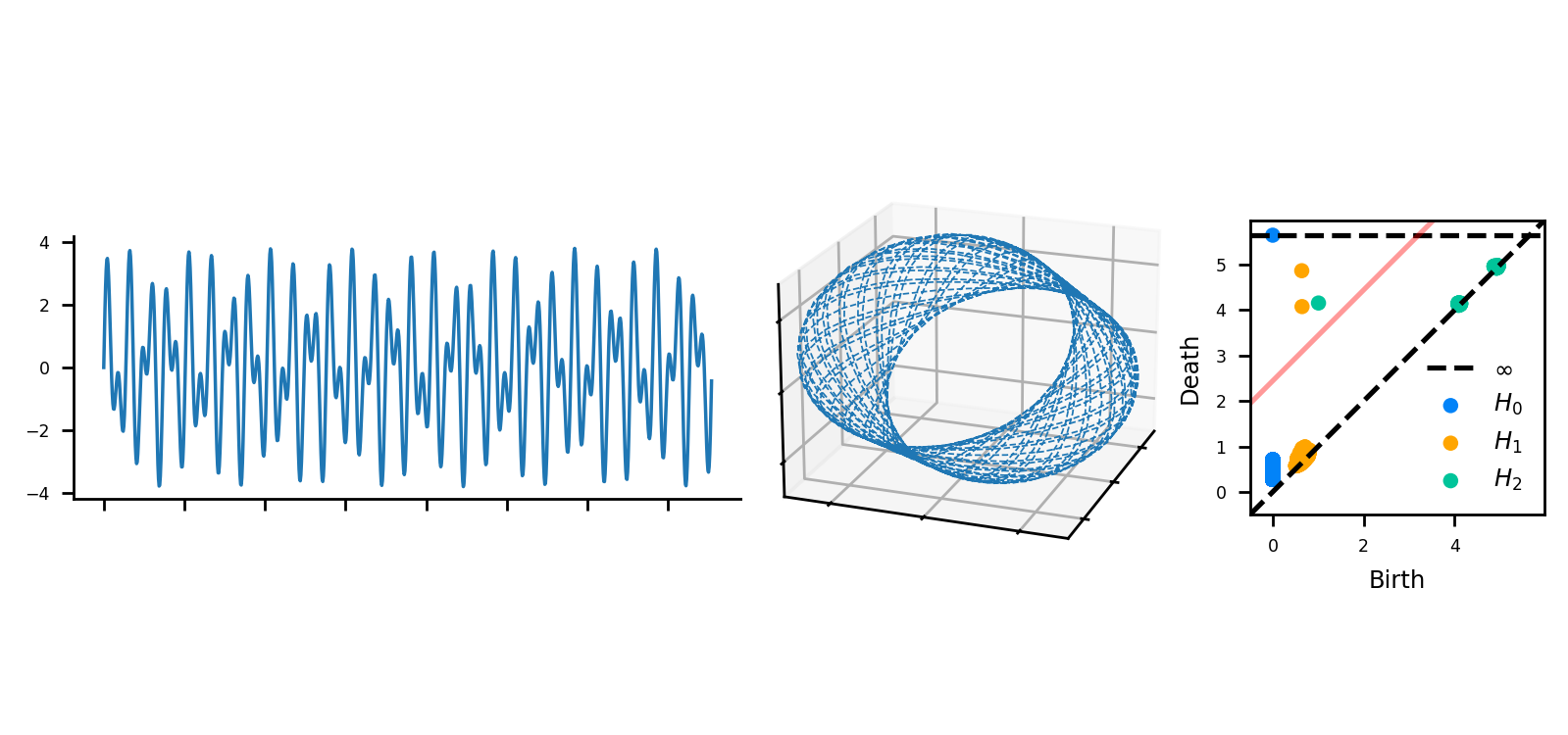}\vspace{-1.2cm}
\includegraphics[width=0.8\textwidth]{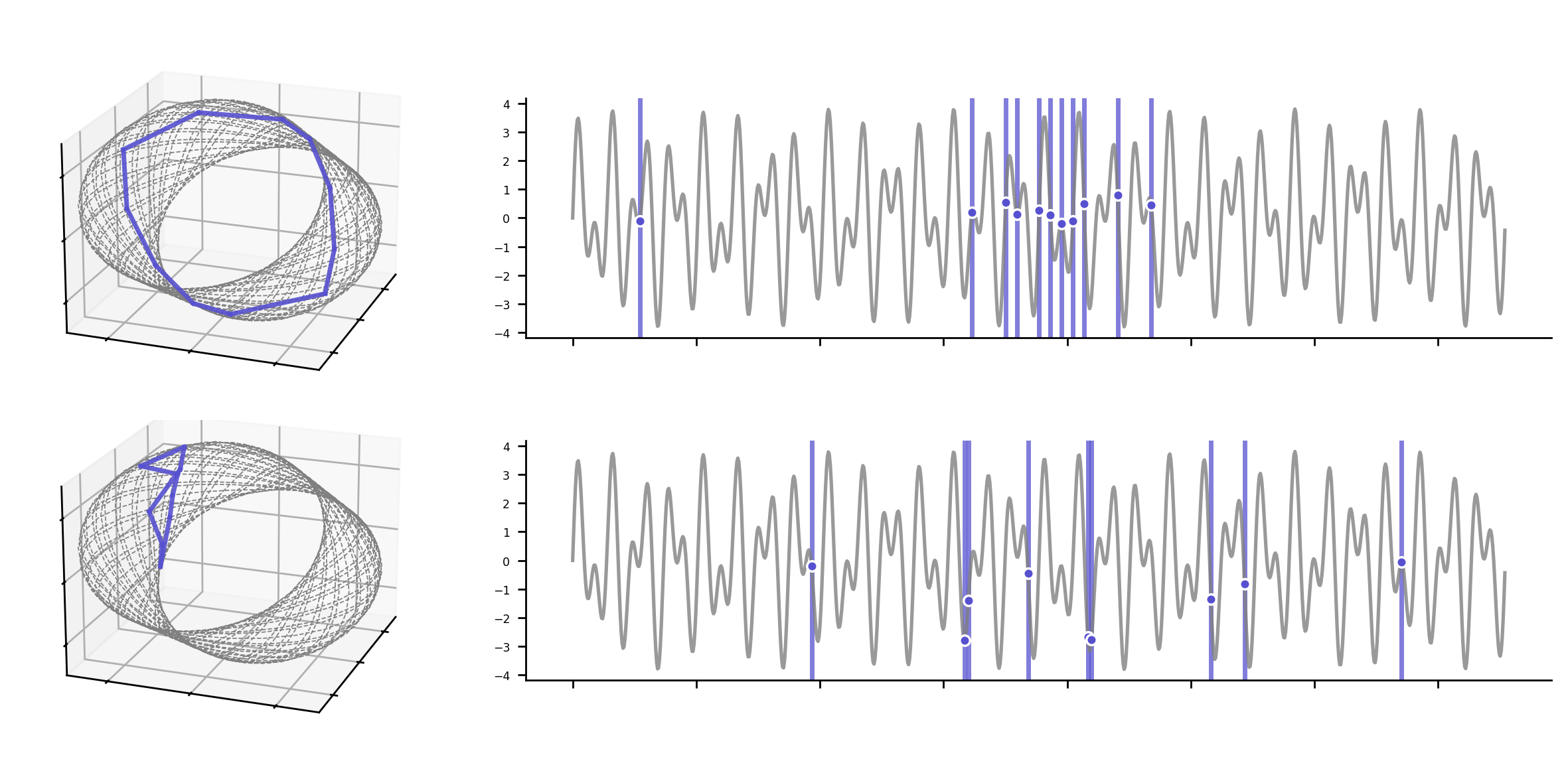}
\includegraphics[width=0.8\textwidth]{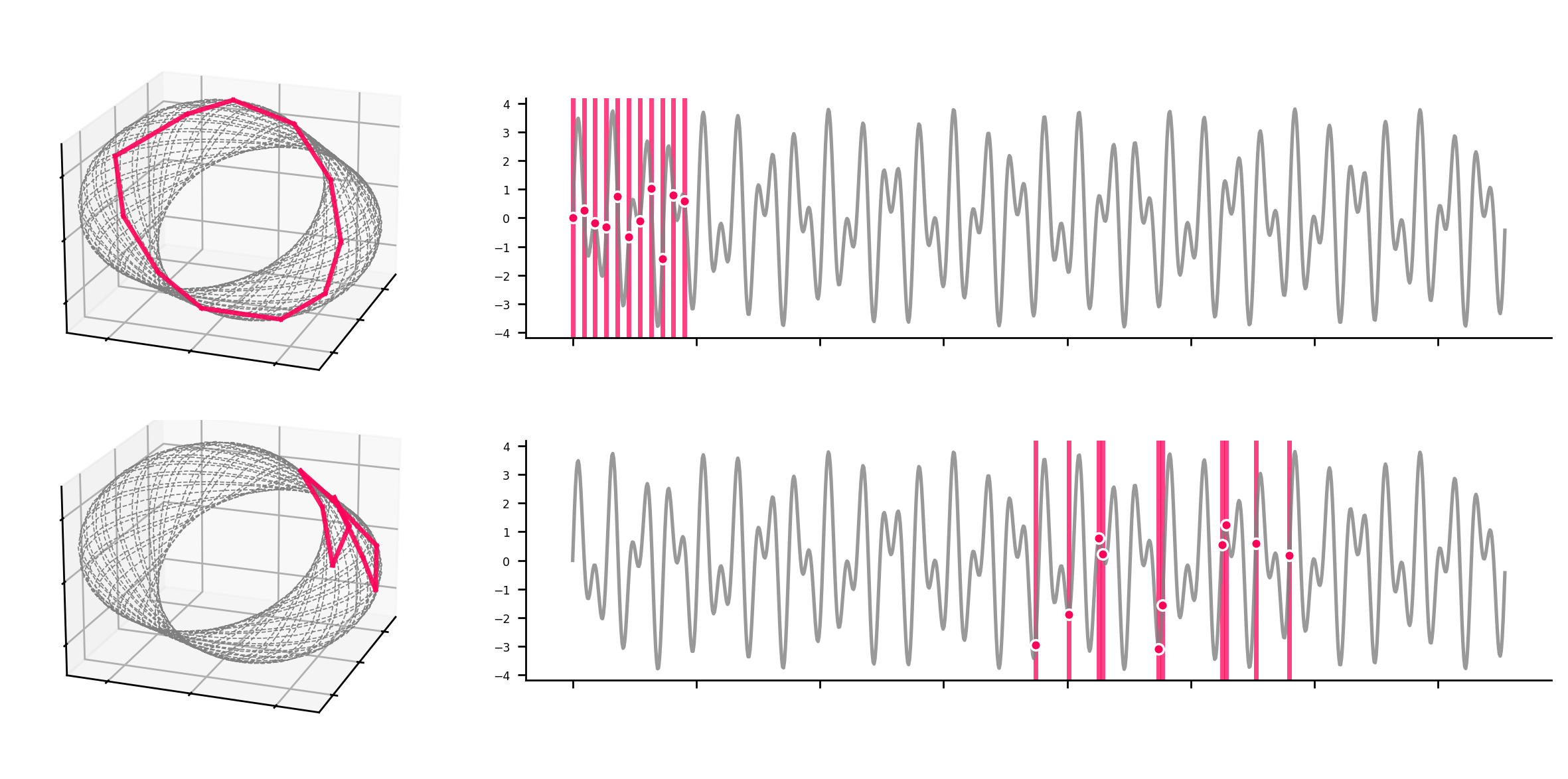}
\caption{We illustrate the two $1$-cycle PH representatives for the optimal embedding of $f(t) = 2\sin{t} + 1.8sin{\sqrt{3}t}$ in the point cloud (left), as well as in the original time series (right).
Top shows the original signal (A) along with the optimal embedding (B) and the corresponding PD (C) with the lower bound.
$2$nd and $3$rd row from top: simplex-based time-optimal cycles. $4$th and $5$th row from top: vertex-based time-optimal cycles.}\label{F:optimal cycles double-sine}
\end{figure}

 \begin{figure}[h!]
 \centering
\includegraphics[width=0.8\textwidth]{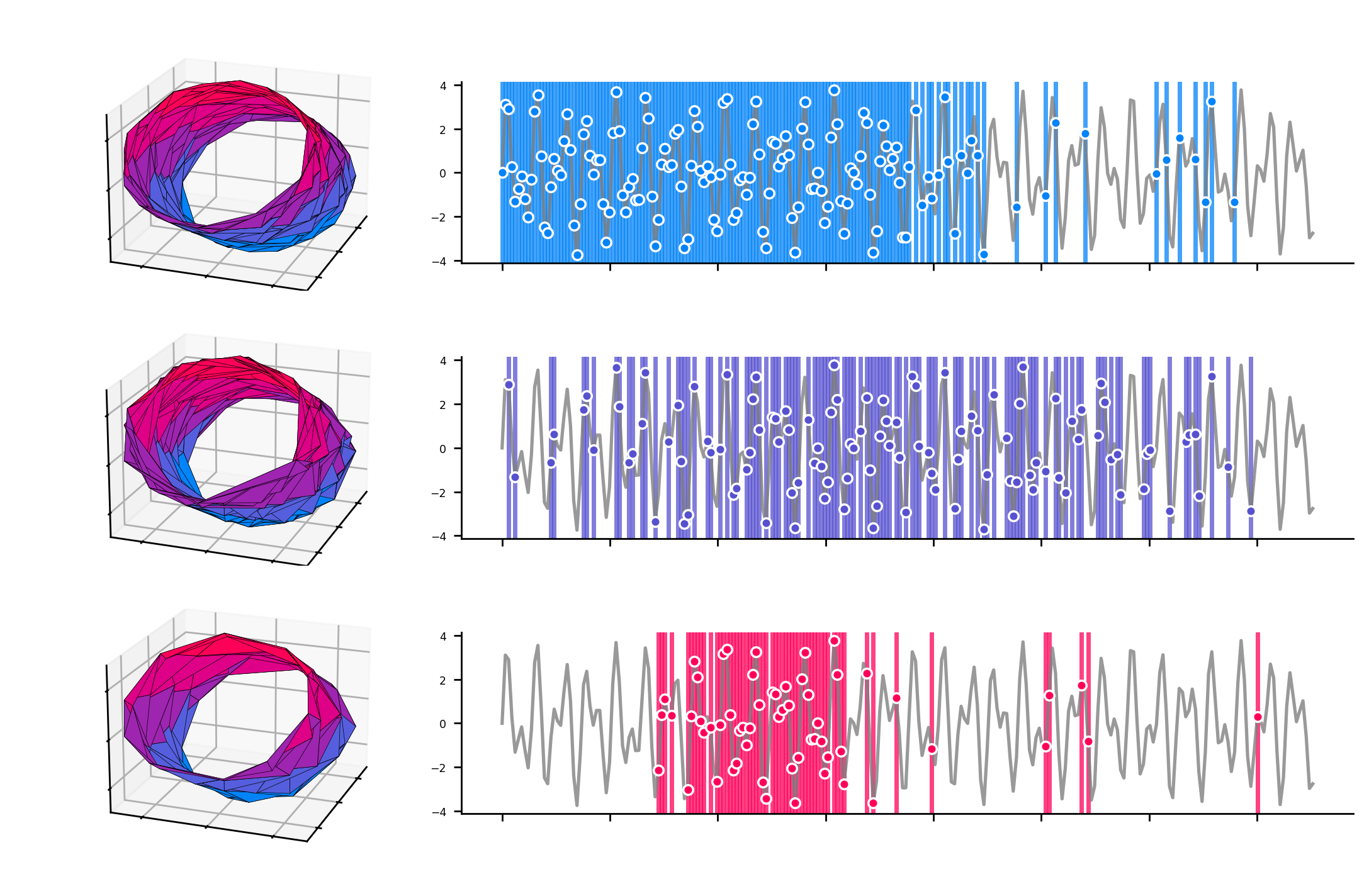}
\caption{Time-optimal representative $2$-cycles for the only significant PH class in $PD_2$, for the optimal embedding point cloud from Figure \ref{fig:example-doublesine-embedding}. Top:  initial representative, given as output by the implementation. Middle: simplex-based optimization. Bottom: vertex-based optimization. We note that the embedding is $4$-dimensional and the triangulation shown on the left column is a  projection onto the first three principal coordinates.}\label{F:optimal cycles double-sine h2}
\end{figure}

We illustrate the results of our computations in Figure \ref{F:optimal cycles double-sine}. 
We observe that both simplex-based and vertex-based representatives demonstrate remarkable consistency in identifying cycle representatives that span approximately one period of the underlying signal.

Our approach can also be used for finding representatives for PH classes in degrees higher than $1$. We illustrate time-optimal PH representatives for the single significant PH class in $PD_2$ in 
Figure \ref{F:optimal cycles double-sine h2}.

We believe that such time-optimal representatives are critical for meaningful signal interpretation, not only for the reasons elucidated earlier, namely the need to be able to physically interpret the results, e.g., in a dynamical systems application, but also because, as the example in Figure \ref{F:optimal cycles double-sine h2} illustrates, visualising cycle representatives in the embedding point clouds is a challenging problem: while visualising appropriately $1$-cycle representatives is challenging for point clouds in more than $3$-dimensions, having any reasonable visualisation in the embedding point cloud becomes an even  more challenging task for cycle representatives in higher homology degrees.

\subsection{Delayed oscillator models of El Nin\~o Southern Oscillation}
We consider quasi-periodic time series arising from a delayed oscillator model of the El Nin\~o Southern Oscillation \cite{CBCP17}. In particular, the model depends on a parameter $\kappa$ encoding the strength of the ocean-atmosphere coupling, and we study three different time series for the values $\kappa\in \{1.4, 1.65, 1.9\}$. We provide details about this data set in Appendix \ref{A:enso}.

\begin{figure}[h!]
\centering
\includegraphics[width=0.8\textwidth]{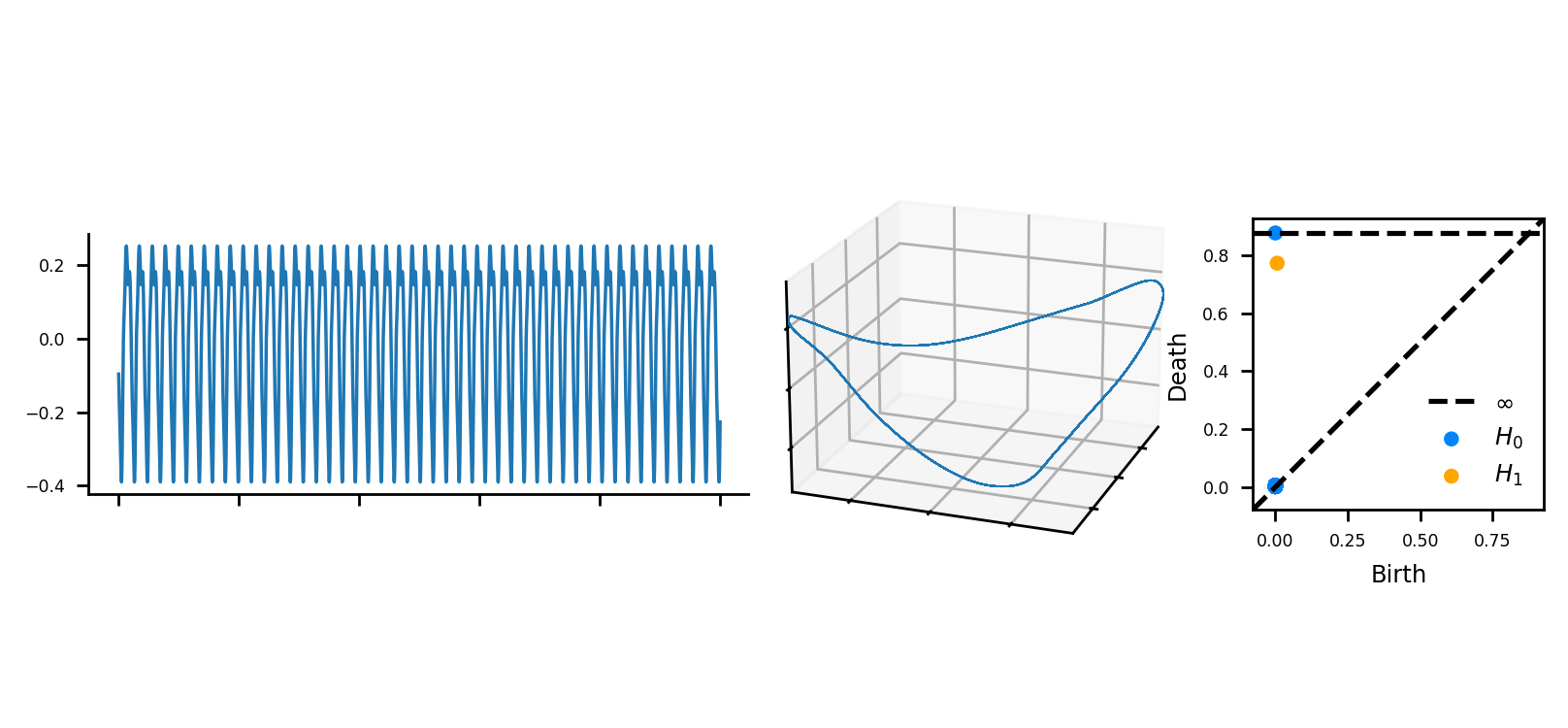}\vspace{-0.5cm}
\includegraphics[width=0.8\textwidth]{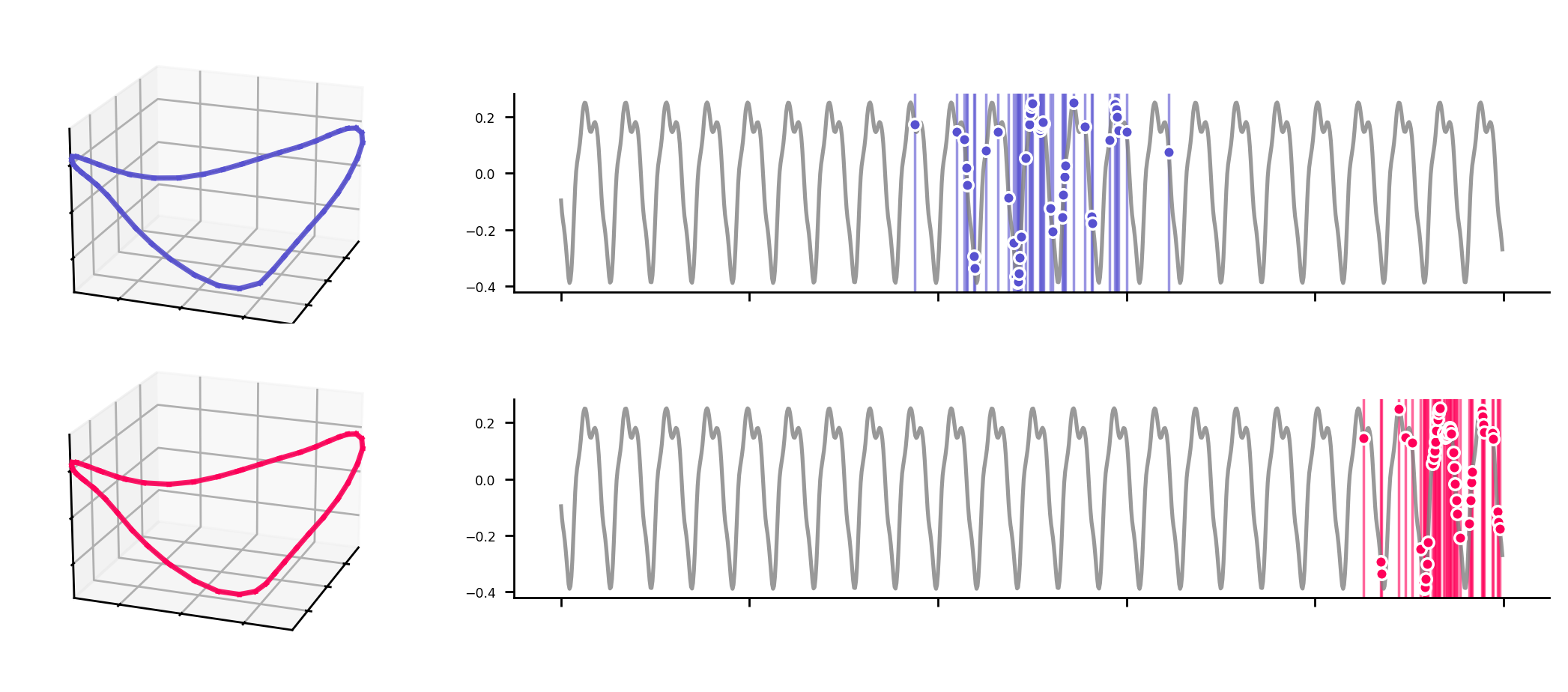}
\caption{Computations for ENSO model for $\kappa=1.65$. Top: original time series, optimal embedding and resulting PD. Middle: simplex-wise time representatives, with relaxed persistence. Bottom:  vertex-wise time representatives, with relaxed persistence. }\label{F:ENSO 1.65}
\end{figure}

\begin{figure}[h!]
\centering
\includegraphics[width=0.8\textwidth]{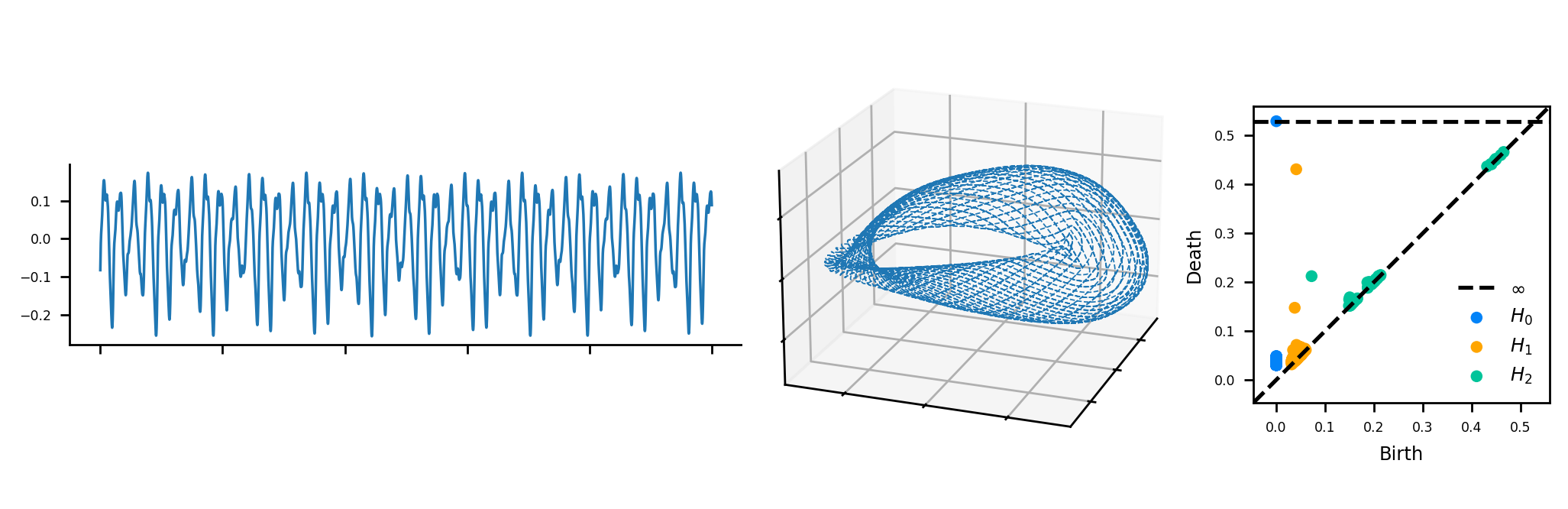}
\includegraphics[width=0.8\textwidth]{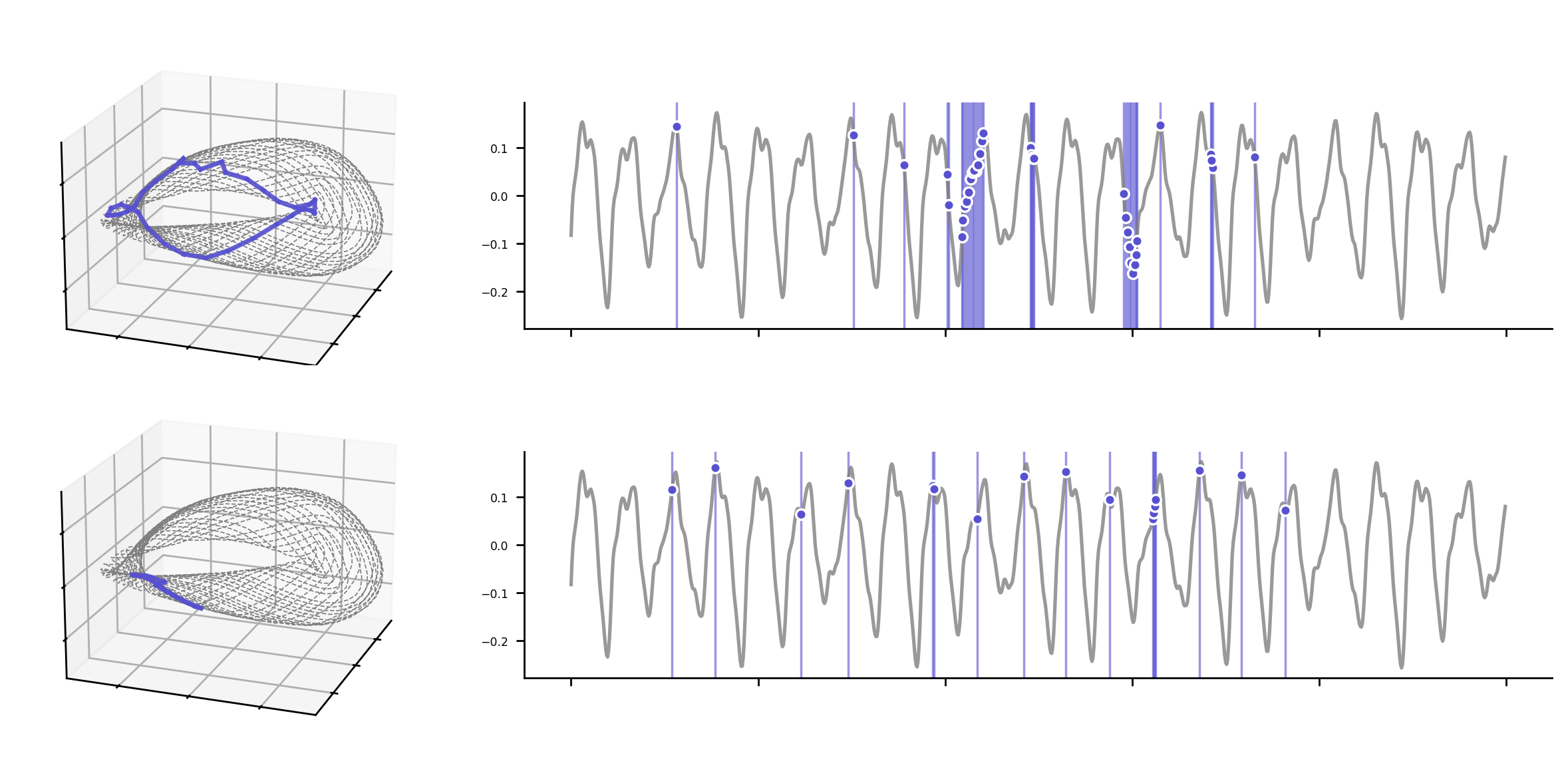}
\includegraphics[width=0.8\textwidth]{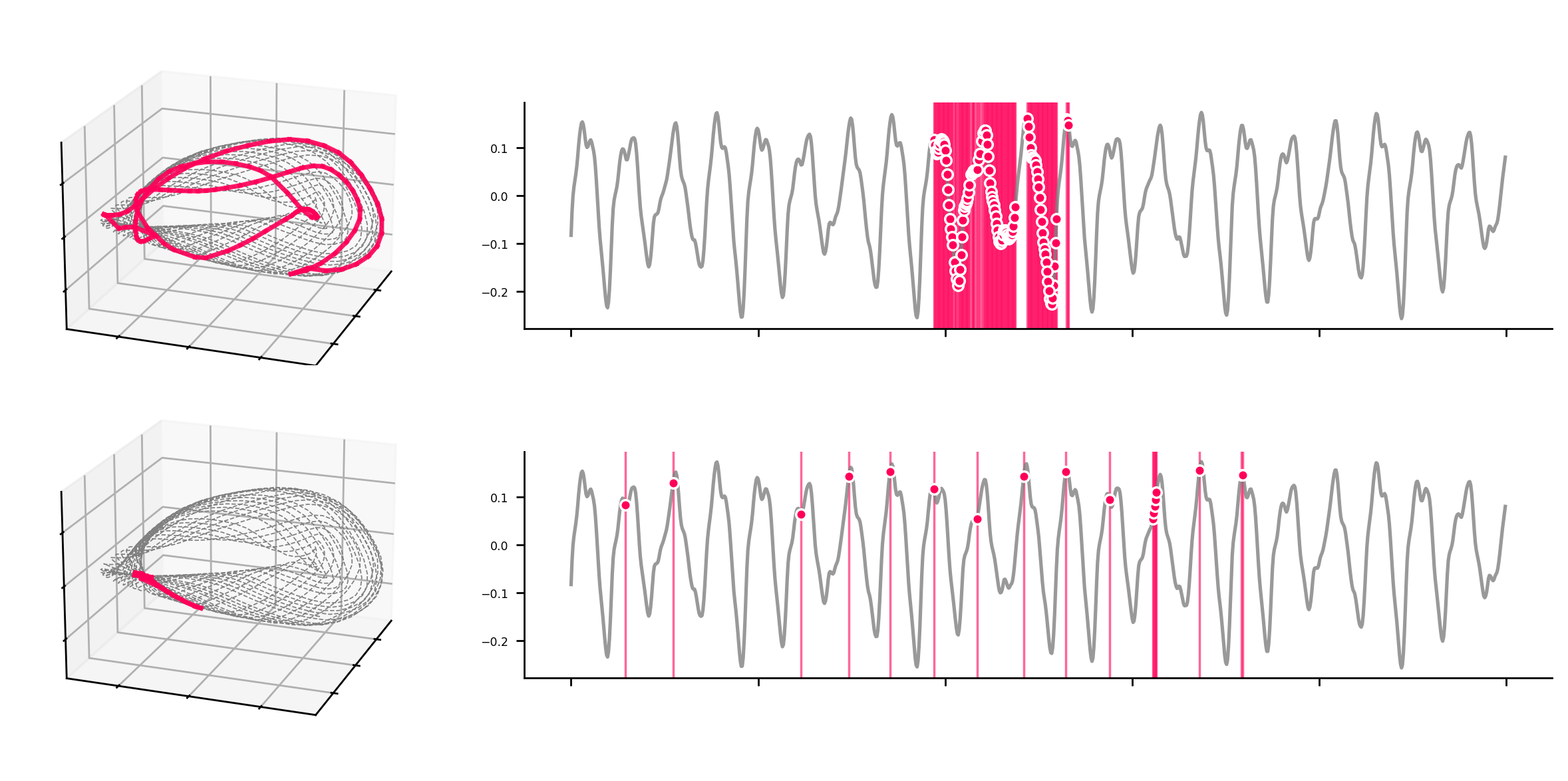}
\caption{Computations for ENSO model for $\kappa=1.4$. Top: original time series, optimal embedding and resulting PD. Middle: simplex-wise time representatives, with relaxed persistence. Bottom:  vertex-wise time representatives, with relaxed persistence. }\label{F:ENSO 1.4}
\end{figure}

\begin{figure}[h!]
\centering
\includegraphics[width=0.8\textwidth]{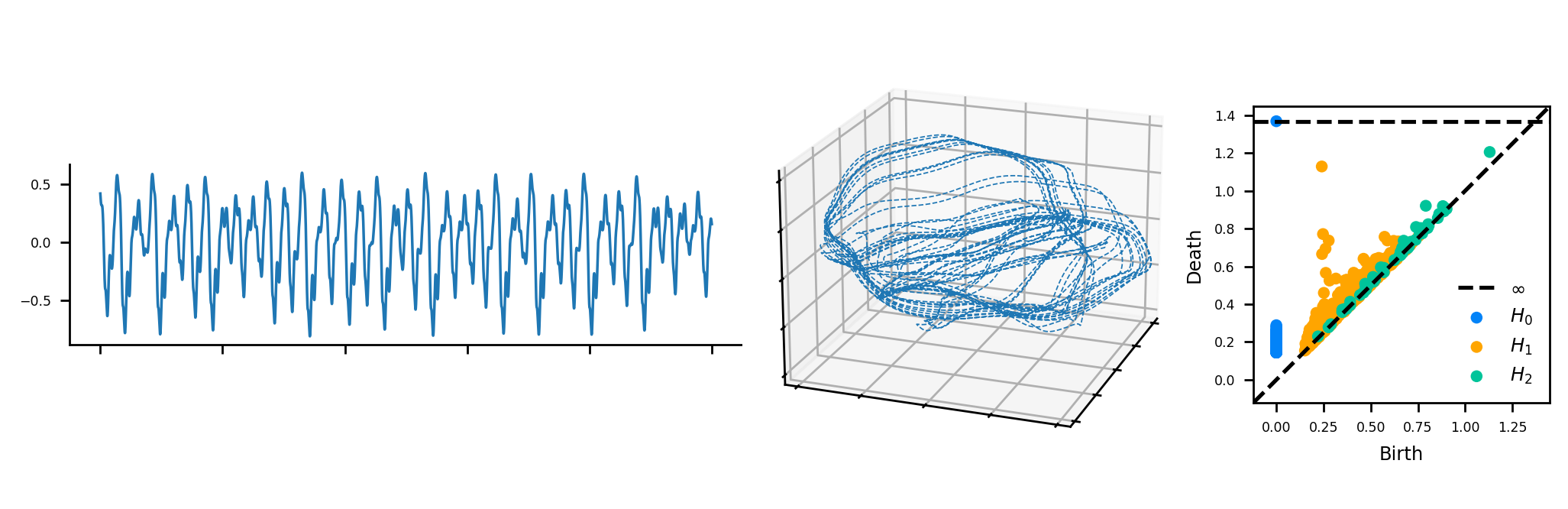}
\includegraphics[width=0.8\textwidth]{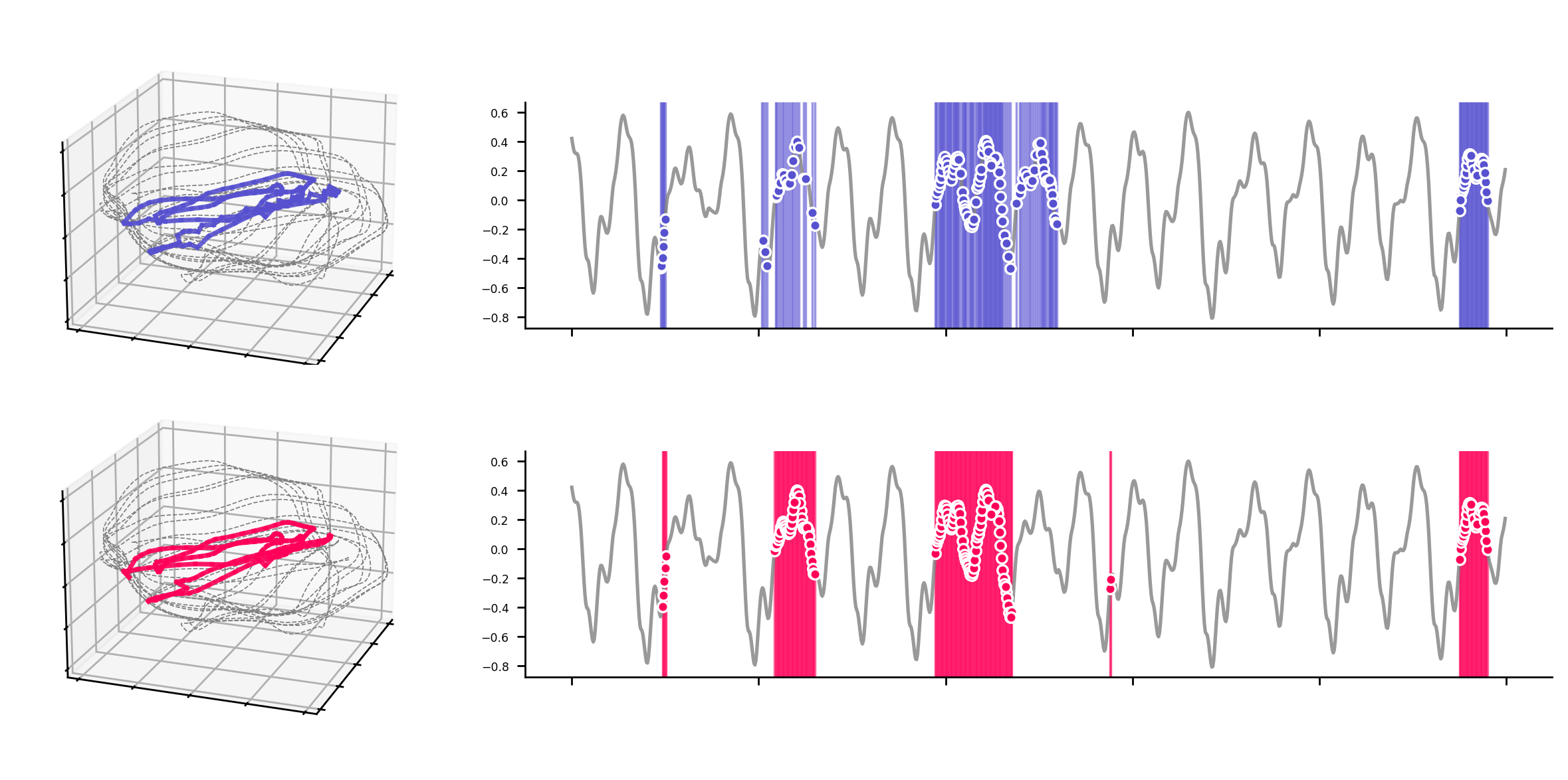}
\caption{Computations for ENSO model for $\kappa=1.9$. Top: original time series, optimal embedding and resulting PD. Middle: simplex-wise time representatives, with relaxed persistence. Bottom:  vertex-wise time representatives, with relaxed persistence. }\label{F:ENSO 1.9}
\end{figure}

Given such a time series, we first compute its embedding parameters by using the methodology developed in \cite{perea2019topological}. We give more details about the parameter computations in Appendix \ref{A:optimal par}. We then compute persistent homology with the Vietoris-Rips complex of the resulting points clouds, and for each point in the persistence diagram $PD_1$, we compute time-optimal PH representatives. The optimization problem is performed with an evenly spaces subsample of 500 points due to computational constraints.  We illustrate the results of the computations in Figures \ref{F:ENSO 1.65}, \ref{F:ENSO 1.4} and \ref{F:ENSO 1.9}.

To be able to obtain representatives that are easier to interpret for our application, we compute a relaxation of time-optimal PH representatives, in which we allow simplices to have later birth times. We provide more details about this relaxation technique in Section \ref{A:approx}. Here we report results for the relaxed version of the PH cycle representatives that have $90\%$ persistence of the corresponding class (as opposed to full persistence). We report the results for cycle representatives with full persistence in Appendix \ref{A: full PH representatives}.

We believe that our approach offers a new perspective for climate science, providing a sophisticated method to characterize oscillatory behaviors that goes beyond the standard analysis of state space embeddings.
More specifically we observe, in certain coupling scenarios, cycle representatives that cluster tightly within specific temporal windows yet have very dispersed spatial configurations (see bottom row of Figure \ref{F:ENSO 1.4}). We aim to further investigate the insights that our methods can bring in  the study  of the El Nin\~o Southern Oscillation in future work.

\section{Conclusions and open questions}\label{S: conclusion}
In the present work we have introduced algorithms for computing PH representatives that correspond to data points that are close in time, in an appropriate sense. We have illustrated the outputs of the different algorithms with synthetic quasi-periodic signals, and several univariate time series arising as models of the El Ni\~no Southern Oscillation.
In most examples that we consider, both methods (simplex and vertex based) manage to extract cycle representatives that are within one period of oscillation, while, overall vertex-based optimization performs better.

One key aspect in extracting representatives that are physically meaningful is the relaxation of birth times for representatives. We aim to further investigate such relaxation techniques, as well as their stability, in future work. In particular, we note that such relaxation techniques are widely applicable to many different notions of optimal PH cycle representatives. 

Our methods can be used to compute time-optimal $p$-cycles, for any $p\geq 0$. In the current work we focused on $p=1$, but  we also provide  an example for $p=2$, for the synthetic quasi-periodic time series, see Figure \ref{F:optimal cycles double-sine h2}. We aim to optimise our implementation in future work to be able to compute higher-degree time-optimal cycles for the ENSO model time series.
Being able to compute such higher-degree PH representatives is particularly important in the study of quasi-periodic time series, which have many topological features of interest in homological degrees greater than $1$.
We also note that  it might be of interest to explore to which extent the choice of norm in the simplex-based method plays a role in the computation of time-optimal cycles. 

Finally, we note that our methods can be applied to other types of time-varying data, including time-varying networks and point clouds.
There are several notions of time-varying networks; we refer the reader to \cite{HOLME201297} for an overview. If one considers time-varying networks in which edges may be deleted or added over time, then one does in general not obtain a nested sequence of simplicial complexes, but rather a zigzag of simplicial complexes, for which one can compute PH with  the zigzag  algorithm \cite{CdSM09,MMK23,KMS20}.  
In particular, we note that for the study of temporal networks, one might be interested in a notion of time dispersion different from the one considered here for time series, and is closer in spirit to the simplex-based approach; namely, one might be interested in defining the time dispersion of a $p$-chain as the difference between minimum and maximum  time labels of any of the $1$-simplices (i.e., edges in the network) contained in the chain.
Zig-zag algorithms have also been used to study time-varying point clouds, for instance in the analysis of flocking behaviour \cite{KMS20}. We will study applications of our methods to such types of data sets in future work.

\section*{Acknowledgments}
We thank Hannah Christensen for having shared with us  the code to compute the ENSO models time series, and for many enlightening discussions.

\bibliography{lipics-v2021-sample-article}

\newpage
\appendix

\section{(Persistent) homology with coefficients over arbitrary fields}\label{A:PH}

The main ingredient in defining simplicial homology and persistent homology over arbitrary coefficient fields is given by a notion of \emph{orientation} of simplices. 

An \define{orientation} of a simplex is the choice of a total order on its vertices. Two orientations of a simplex are defined to be equivalent if they differ by an even permutation.

\begin{figure}[h!]
\begin{tikzpicture}
\node (a)at (0,0) {$a$};
\node (b)at (1,0) {$b$};
\node (c)at (0.5,1) {$c$};
\draw[->] (c)edge(a)
(b) edge  (a)
(c)edge(b);
\end{tikzpicture}\qquad \qquad 
\begin{tikzpicture}
\node (a)at (0,0) {$a$};
\node (b)at (1,0) {$b$};
\node (c)at (0.5,1) {$c$};
\draw[->] (a)edge(c)
(a) edge  (b)
(c)edge(b);
\end{tikzpicture}
\qquad \qquad 
\begin{tikzpicture}
\node (a)at (0,0) {$a$};
\node (b)at (1,0) {$b$};
\node (c)at (0.5,1) {$c$};
\draw[->] (a)edge(c)
(a)edge(b) 
(b)edge(c);
\end{tikzpicture}
\caption{A simplex with different orientations, where given an order $<$ on the vertices, we denote $i<j$ by drawing an arrow $j\to i$. Left: orientation given by lexicographical order $a<b<c$. Middle: orientation given by $b<c<a$. Right: Orientation given by $c<b<a$. This orientation is equivalent to the one given by the lexicographical order.}\label{F: orientation simplex}
\end{figure}
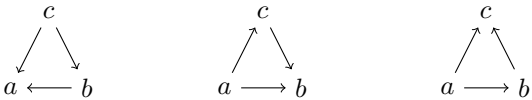

In practice, one usually chooses a total order on the set of vertices of a simplicial complexes, and then gives to each simplex the orientation induced by the order on its vertices. We illustrate an example of a simplex with different orientations in Figure \ref{F: orientation simplex}.

The main difference to defining simplicial homology over the field with $2$ elements then relies on how the boundary maps on chains are defined. One defines
\begin{alignat}{2}
\bd_p: C_p(K) &\notag\longrightarrow C_{p-1}(K)\\
\sigma=(x_0,\dots , x_p) &\notag \mapsto \sum_{i=0}^p (-1)^i  (x_0,\dots , \widehat{x}_i, \dots, x_p)
\end{alignat}
where the symbol $\widehat{x}_i$ means that the $i$th vertex has been deleted.

We refer the reader to \cite{hatcher2002algebraic} for more details on simplicial homology over general coefficient fields.

\section{Details on algorithms and implementation}\label{A:details algo impl}
\subsection{Optimal homologous cycle}
Given an initial p-cycle $c_0 \in C_p(K)$, finding an homologous cycle that minimizes a given loss function $\ell$ is defined as:
\begin{align*}
\text{min}   & \quad \ell(c) \\
\text{subject to} & \quad c = c_0 + \bd_{p+1}(w) \\
                  & \quad w \in C_{p+1}(K)
\end{align*}
The search space consists of the initial cycle plus the boundary of a higher dimensional simplex. Guaranteeing, by definition, that the solution $c$ is homologous to the initial cycle $c_0$.

In a simplicial filtration, finding a representative for a homology class with birth and death $(b,d)$ requires that the optimal representative also shares the same birth and death.
This reflects in restricting the search to the simplicial complex existing at birth time ($K_b$).
Which resolves to the following problem \cite{CF08}:
\begin{align*}
\text{min}   & \quad \ell(c) \\
\text{subject to} & \quad c = c_0 + \bd_{p+1}(w) \\
& \quad w \in C_{p+1}(K_b) \\
& \quad \bd_{p+1}(w) \in C_p(K_b)
\end{align*}

Where the search is contained in cycles composed of simplices that have a birth \textit{before or at} the birth of the homology class.
The third constraint is redundant since if $w \in K_b$ then $\bd(w) \in K_b$, that is, the birth of a simplex is always later than or equal to the ones that make up its boundary. 

\subsection{Linear Programming}\label{A:lin progr}
If the function we optimize is linear, then it is possible to solve the previous problem with linear optimization \cite{escolar2016optimal}.
The added step is separating the coefficients vector $c$ into a positive $c^+$ and a negative $c^-$ part such that $c = c^+ - c^-$ allowing us to set both $c^+,c^- \geq 0$.

The optimization can be done both by allowing only binary coefficients $(\mathbb{F}_2)$ which requires a solver capable of mixed linear programming or can be relaxed to $\mathbb{R}$.
In general using the $l_1$ norm over $\mathbb{R}$ is sufficient to provide sparse solutions \cite{lu2021minimal}.

Let $\mathbf{c} =\mathbf{c}^+ - \mathbf{c}^- $ be the vector of coefficients in $\mathbb{R}$.
Let $W$ be a non-negative weight matrix, then the linear programming formulation is the following:
\begin{align*}
\text{min}   &  \quad \lVert  W \mathbf{c} \rVert_1 =  \sum_i \sum_j w_{ij} (c_{j}^+ + c_{j}^-) \\
\text{subject to  } & (\mathbf{c}^+ - \mathbf{c}^-) = \mathbf{c_0} + \bd_{p+1}(\mathbf{w}) \\
            & \mathbf{w}\in \mathbb{R}^m\\
            & \mathbf{c} \in \mathbb{R}^n\\
            & \mathbf{c}^+,\mathbf{c}^- \geq 0
\end{align*}

The initial representative $\mathbf{c}_0$ can be obtained by most persistent homology methods.
It appears naturally from the decomposition $R= \bd_{k} V$ of the filtration boundary matrix.
Given the birth $b$ of the homology class we consider the sets of $p$ and $p+1$ simplices that are alive at filtration time $b$:
$$
P = \{\sigma \in S_{p}(K_b) \mid \text{birth}(\sigma) \leq b\}\\
$$
$$
Q = \{\sigma \in S_{p+1}(K_b) \mid \text{birth}(\sigma) \leq b\}\\
$$
The solution of the previous problem with the rows and columns of the boundary operator $\bd_{p+1}$ restricted to $\bd_{p+1}[P,Q]$ is a cycle with persistence $(b,d)$.

A substantial speedup comes from restricting the domain of the boundary operator even further by considering the set \cite{lu2021minimal}:
$$
\hat{Q} = \{\sigma  \mid R[:,\sigma] \neq 0\}  
$$
Which are the nonzero columns of the $R$ matrix resulting from the previous decomposition $R = \bd_{p+1}V$.
This significantly reduces the number of conditions without affecting the search space.
The resulting problem is the following:
\begin{align*}
\text{min}   &  \quad \lVert  W \mathbf{c} \rVert_1 = \sum_i \sum_j w_{ij} (c_{j}^+ + c_{j}^-) \\
\text{subject to  } & (\mathbf{c}^+ - \mathbf{c}^-) = \mathbf{c_0} + \bd_{p+1}[P,\hat{Q}](\mathbf{w}) \\
            & \mathbf{w}\in \mathbb{R}^{|\hat{Q}|}\\
            & \mathbf{c} \in \mathbb{R}^{|P|}\\
            & \mathbf{c}^+,\mathbf{c}^- \geq 0 \, .
\end{align*}

\section{Time-delay embeddings of quasi-periodic time series}\label{A:optimal par}

In this section we provide details on a framework  developed in \cite{GP24}, which provides a rigorous methodology to compute  optimal parameters for sliding window embeddings of quasi-periodic functions.  
This computational methodology allows to obtain embedding parameters so that the persistent homology, computed with respect to the Vietoris-Rips complex of the obtained embedding point cloud, reflects in a robust manner the persistent homology of a hypertorus in as many dimensions as the frequencies characterising the quasi-periodic function. We emphasise that there are no contribution by the authors in this section, and all contributions are due to the authors of \cite{GP24}.

\subsection{Quasiperiodic Functions}

\begin{definition}[Quasi-periodic function] Let $N\in \mathbb{N}$ and $\mathbb{T}^N =\left (\mathbb{R}/2\pi \mathbb{Z}\right)^N$.
A function \( f \colon \mathbb{R} \to \mathbb{C} \) is \define{quasi-periodic} if there exists a vector \( \omega = (\omega_1, \dots, \omega_N)\in \mathbb{R}_{\geq 0}^\mathbb{N} \) with components \(\omega_i\) linearly independent over \(\mathbb{Q}\) and a  function \( F : \mathbb{T}^N \to \mathbb{C} \) such that:
$$
f(t) = F(\omega_1 t, \dots, \omega_N t)\, .
$$
The vector $\omega$ is called \define{frequency vector} of $f$, and the function $F$ \define{parent function} of $f$.
\end{definition}
In other words, the function $f$ is of the form
$$
f(t) = c_1e^{i\omega_1t} + c_2e^{i\omega_2t} + \dots + c_Ne^{i\omega_Nt}\, 
$$ 
where $\omega_1, \dots \omega_N$ are linearly independent over $\mathbb{Q}$ (incommensurable) and $c_i>0$.

The sliding window embedding of \(f(t)\) with delay \(\tau > 0\) and embedding dimension \(d+1\) is defined for any $t\in \mathbb{R}$ as:
\[
SW_{d,\tau}f(t) = 
\begin{bmatrix}
f(t) \\ f(t+\tau) \\ \vdots \\ f(t+d\tau)
\end{bmatrix} \in \mathbb{C}^{d+1}.
\]
The geometry of the embedding \(SW_{d,\tau}f\) depends on the parameters \(\tau\) and \(d\), as well as the properties of \(f\).

For practical computations, \(f(t)\) is approximated by its truncated Fourier series. Given positive integers $N$ and $K$, one considers  a restricted integral square box of side length bounded by $2K$ in an $N$-dimensional grid with integer grid points:
$$
I^N_K = \{\mathbf{k} \in \mathbb{Z}^N \mid \lVert \mathbf{k} \rVert_\infty < K \} \, .
$$
One then  defines the truncation:
$$
S_K f(t) = \sum_{\mathbf{k} \in I^N_K} \hat{F}(\mathbf{k}) e^{i\langle \mathbf{k}, \omega t \rangle},
$$
where the integer $K$ can be thought of as controlling the fidelity of the truncation.

\subsection{Optimal dimension}
We can estimate the nonzero Fourier coefficients $\hat{F}(\mathbf{k})$ and their frequency locations $\langle \mathbf{k}, \omega \rangle$  as follows:
$$
\text{supp}(\hat{F}_K) := \left\{ \mathbf{k} \in \mathbb{Z}^N \ \middle|\ \hat{F}(\mathbf{k}) \neq 0 \text{ and } \|\mathbf{k}\|_\infty \leq K \right\}\, .
$$
We take the embedding dimension $d$ to be the cardinality of supp$(\hat{F}_K)$, (the number of prominent peaks in the spectrum of $f$).

\subsection{Optimal delay}
The choice of delay $\tau$ influences the appearances of the homological features in given degrees of the hypertorus in the sliding window embedding. In particular, poor choices can obscure them. Figure \ref{fig:example-doublesine-embedding} shows an example of how a suboptimal delay parameter (bottom) squashes the homology groups.
The sliding window embedding, for dimension $d$ and delay $\tau$, of the truncated Fourier approximation of $f$ is the following
    \begin{align}
        SW_{d,\tau}S_Kf(t) &=\notag
        \begin{bmatrix}
            1 & \dots  & 1 \\
            e^{i \langle \textbf{k}_1,\omega \rangle \tau} &  \dots  & e^{i \langle \textbf{k}_{\alpha},\omega \rangle \tau} \\
            \vdots & \vdots  & \vdots \\
            e^{i \langle \textbf{k}_1,\omega \rangle \tau d} &  \dots  & e^{i \langle \textbf{k}_{\alpha},\omega \rangle \tau d}
        \end{bmatrix}\cdot
        \begin{bmatrix}
            \hat{F}(\textbf{k}_1)e^{i\langle \textbf{k}_1,\omega \rangle t} \\
            \vdots\\
            \hat{F}(\textbf{k}_{\alpha})e^{i\langle \textbf{k}_{\alpha},\omega \rangle t}
        \end{bmatrix}\\
       &\notag  = \Omega_{K,f} \cdot x_{K,f}(t)
    \end{align}

Note that only the $\Omega_{K,x}$ matrix depends on the delay value $\tau$.
The optimal choice of delay is the one that best improves the conditioning number of $\Omega_{K,x}$ so that no toroidal features are squashed.
This is done by minimizing a scalar function that measures the extent to which columns in $\Omega_{K,f}$ are pairwise orthogonal.

\subsection{Persistence Significance Bounds}\label{A:persistence bounds}
The embedding quality depends on the Fourier approximation (choice of $K$) and the embedding parameters $(d,\tau)$ as well as the smallest singular value  \(\sigma_{\min}\) of \(\Omega_{K,f}\). 
We have that  \cite[Theorem 6.8]{GP24} provides lower bounds on the lifespan $b-a$ of points $(a,b)$ in a persistence diagram. Such bounds can be interpreted as  giving a separation between noise (points with lifespan smaller than the bounds) and signal (points with lifespan greater or equal to the bounds).

\newpage

\section{ENSO model time series}\label{A:enso}
The El Niño–Southern Oscillation (ENSO) is a set of coupled ocean-atmosphere phenomena characterized by an irregular cycle of warming (El Niño) and cooling (La Niña) in the eastern tropical Pacific
along with a corresponding variation in sea level pressure \cite{CBCP17,munnich1991study}.
ENSO significantly impacts global weather patterns and is associated with heavy rain in Peru, drought in Indonesia, intensity of the Indian monsoon and the number of hurricanes in North America \cite{CBCP17}. 
The parameter $\kappa$ represents the coupling strength of the ocean-atmosphere coupling, different values of this parameter simulate a variety of ENSO behaviors.

The ENSO data consists in 3 timeseries, one for each value of $\kappa=\{1.4,1.65,19\}$, each with $100 000$ points representing the modeled oscillation.
Each time series is embedded in $\mathbb{R}^n$ using the optimal parameters.
Below we show the different series and persistence diagrams of optimal (labeled MIN in the figures) and suboptimal embeddings.
The persistent homology is calculated using Ripser \cite{bauer2021ripser} using a subsample of 1000 points \cite{cavanna2015geometricperspectivesparsefiltrations}.

\subsection{ENSO model ($\kappa=1.4$)}

\begin{figure}[!htb]
     \centering
     \begin{subfigure}[b]{0.6\textwidth}
         \centering
         \includegraphics[width=\textwidth]{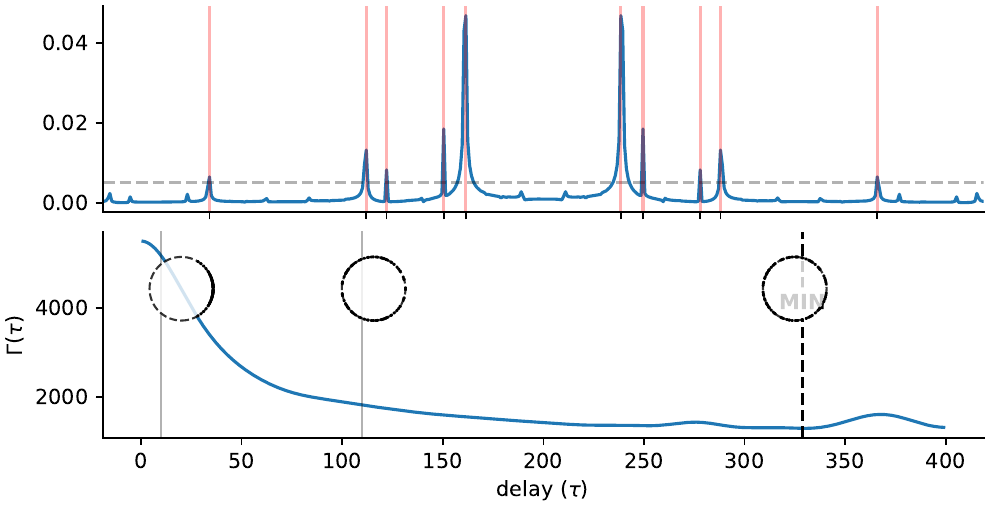}
     \end{subfigure}\vspace{-1.7cm}
     \begin{subfigure}[b]{0.7\textwidth}
         \centering
         \includegraphics[width=\textwidth]{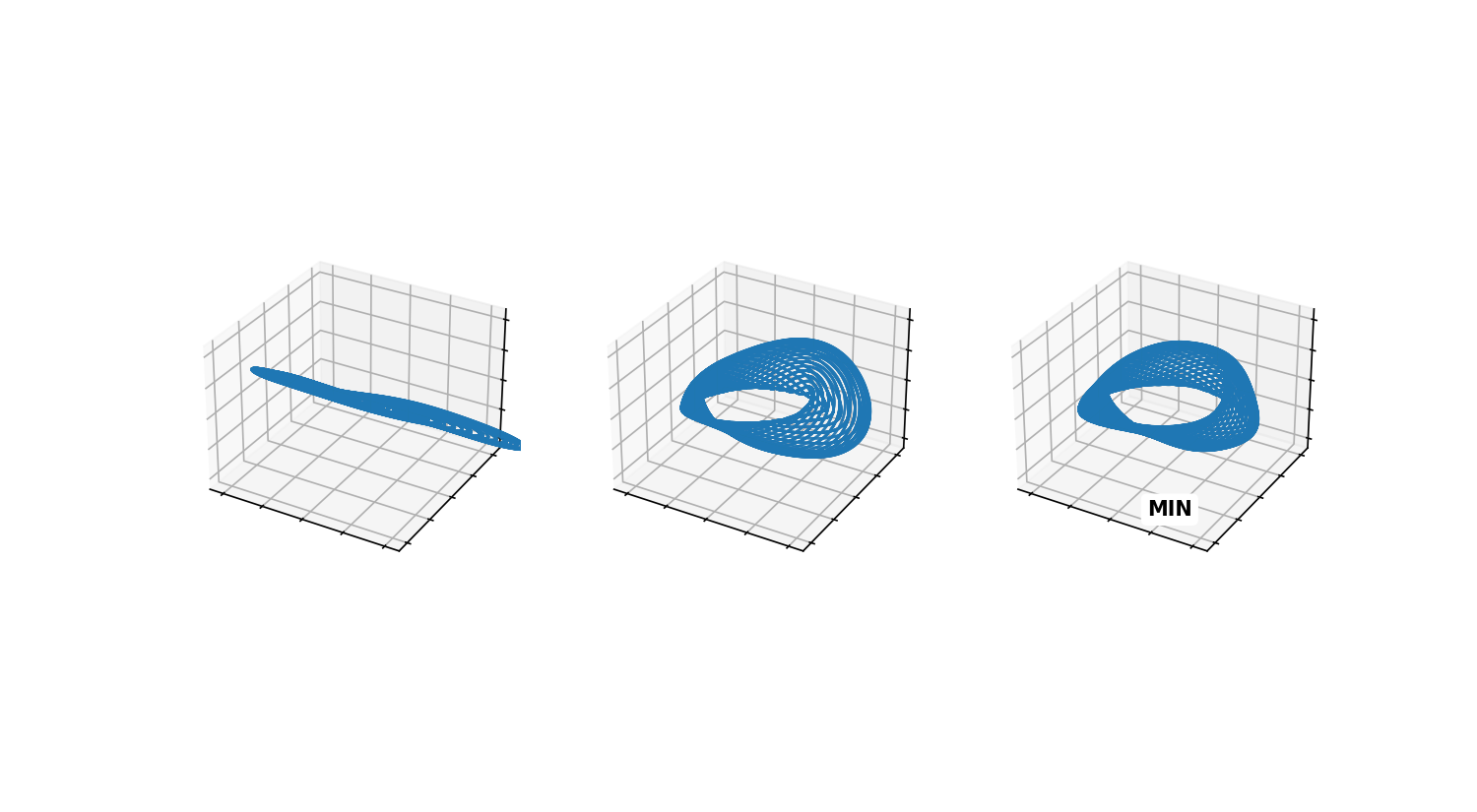}
         \vspace{-3\baselineskip}
     \end{subfigure}\vspace{-0.5cm}
     \begin{subfigure}[b]{0.6\textwidth}
         \centering
         \includegraphics[width=\textwidth]{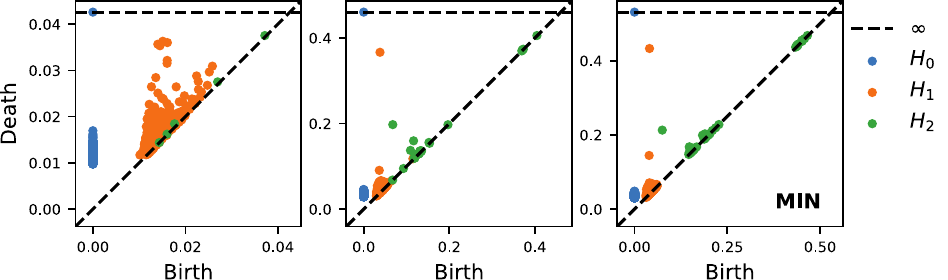}
     \end{subfigure}
        \caption{Optimal sliding window embedding of the ENSO model for the parameter $k=1.4$. \textbf{Top:} Fourier decomposition of the signa, and the average orthogonality of the $\Omega_{K,f}$ matrix (Appendix \ref{A:optimal par}) with respect to a chosen delay value. \textbf{Middle:} PCA reduction of the original embeddings in $\mathbb{R}^{10}$ for optimal (right) and suboptimal delay values (right and middle). \textbf{Bottom:} Persistence Diagram for each embedding.}
        \label{fig:enso_14}
\end{figure}

\newpage
\subsection{ENSO model ($\kappa=1.65$)}
\begin{figure}[!htb]
     \centering
     \begin{subfigure}[b]{0.7\textwidth}
         \centering
         \includegraphics[width=\textwidth]{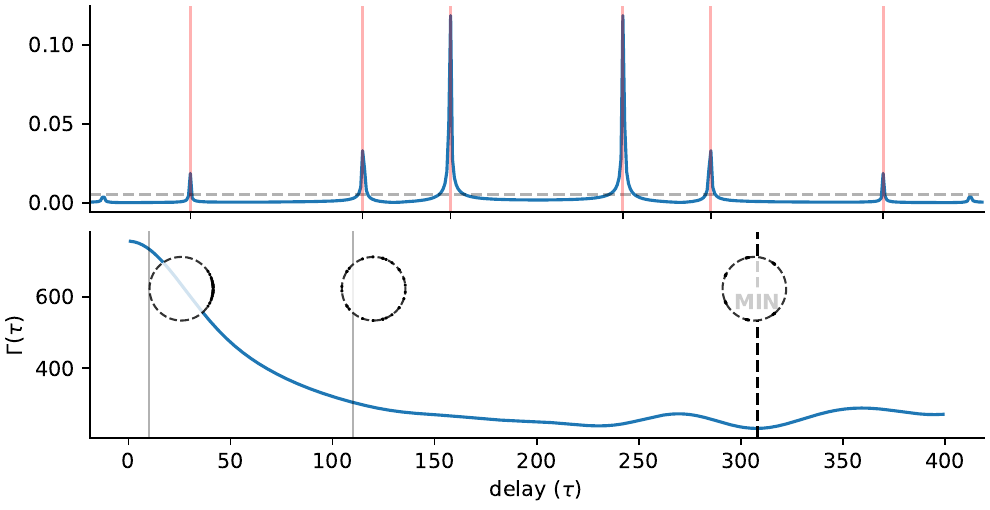}
     \end{subfigure}\vspace{-1.7cm}
     \begin{subfigure}[b]{0.7\textwidth}
         \centering
         \includegraphics[width=\textwidth]{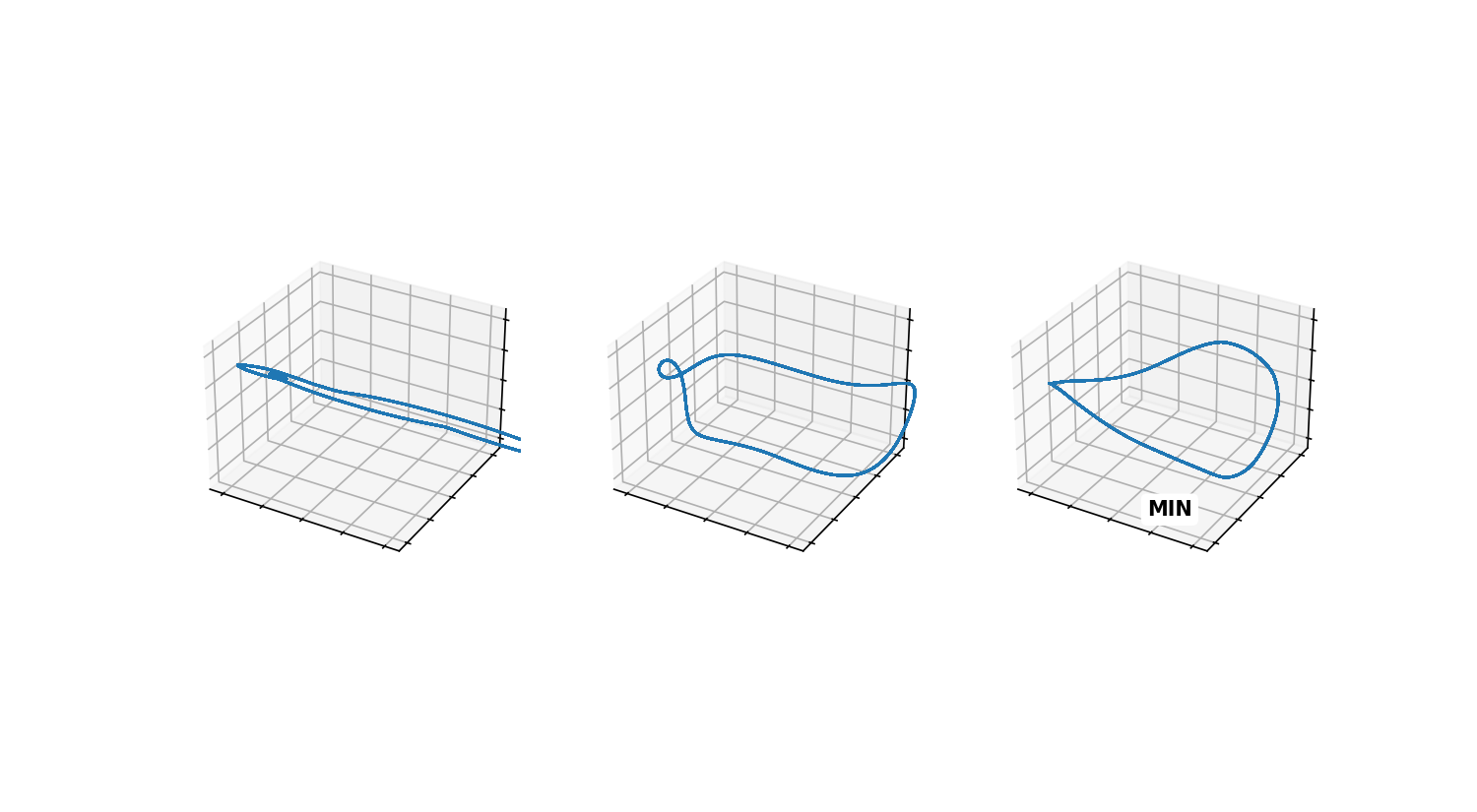}
         \vspace{-3\baselineskip}
     \end{subfigure}\vspace{-0.5cm}
     \begin{subfigure}[b]{0.6\textwidth}
         \centering
         \includegraphics[width=\textwidth]{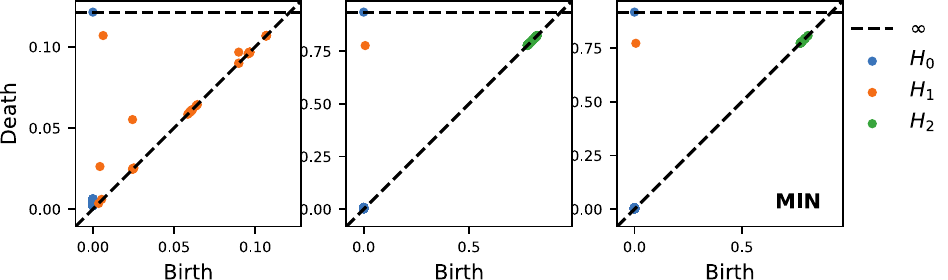}
     \end{subfigure}
         \caption{Optimal sliding window embedding of the ENSO model for the parameter $k=1.65$. \textbf{Top:} Fourier decomposition of the signal, and the average orthogonality of the $\Omega_{K,f}$ matrix (Appendix \ref{A:optimal par}) with respect to a chosen delay value. \textbf{Middle:} PCA reduction of the original embeddings in $\mathbb{R}^{6}$ for optimal (right) and suboptimal delay values (right and middle). \textbf{Bottom:} Persistence Diagram for each embedding.}
        \label{fig:enso_1.65}
\end{figure}

\newpage
\subsection{ENSO model ($\kappa=1.9$)}
\begin{figure}[!htb]
     \centering
     \begin{subfigure}[b]{0.7\textwidth}
         \centering
         \includegraphics[width=\textwidth]{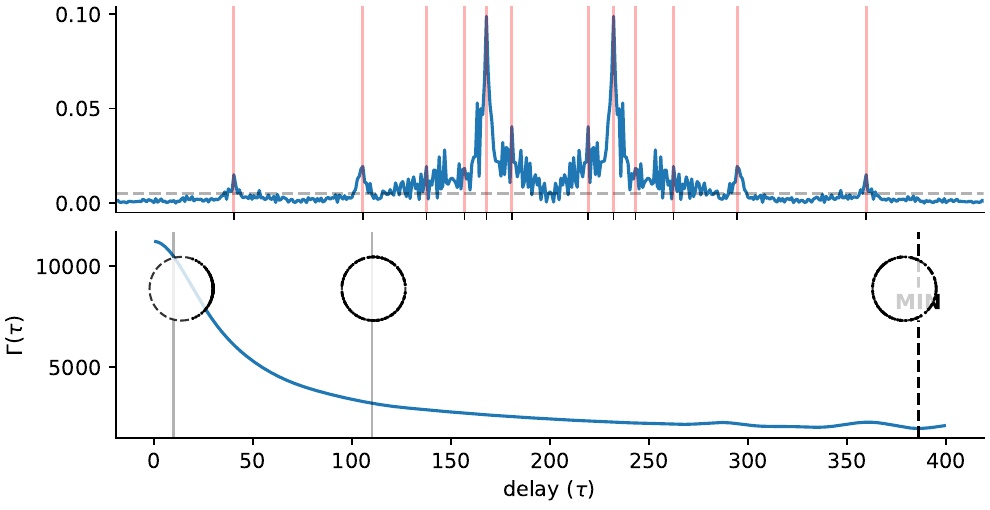}
     \end{subfigure}\vspace{-1.7cm}
     \begin{subfigure}[b]{0.7\textwidth}
         \centering
         \includegraphics[width=\textwidth]{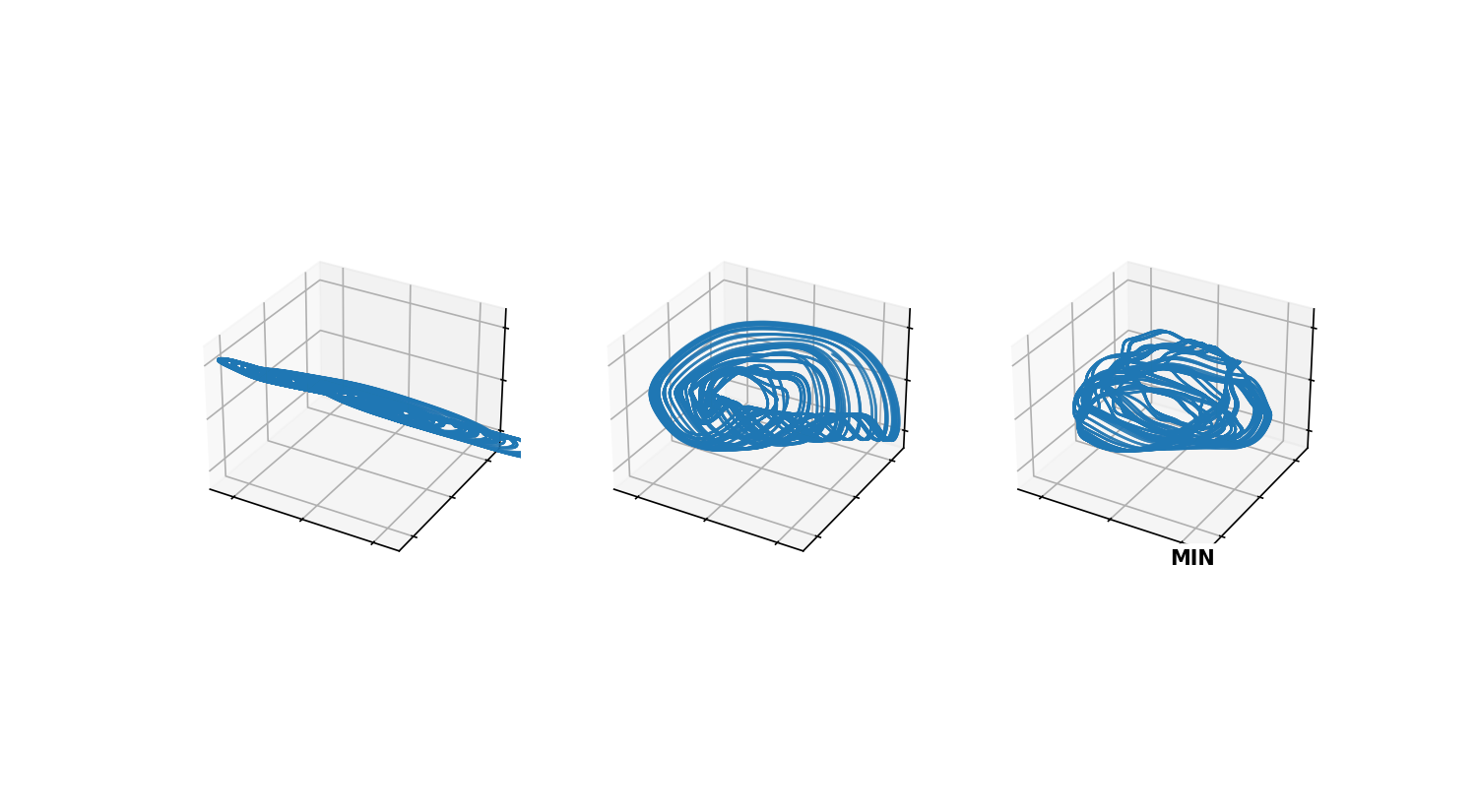}
         \vspace{-3\baselineskip}
     \end{subfigure}\vspace{-0.5cm}
     \begin{subfigure}[b]{0.6\textwidth}
         \centering
         \includegraphics[width=\textwidth]{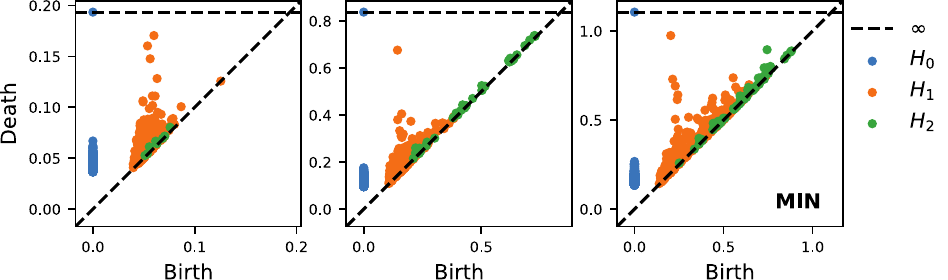}
     \end{subfigure}
        \caption{Optimal sliding window embedding of the ENSO model for the parameter $k=1.9$.\textbf{Top:} Fourier decomposition of the signal, and the average orthogonality of the $\Omega_{K,f}$ matrix (Appendix \ref{A:optimal par}) with respect to a chosen delay value. \textbf{Middle:} PCA reduction of the original embeddings in $\mathbb{R}^{12}$ for optimal (right) and suboptimal delay values (right and middle). \textbf{Bottom:} Persistence Diagram for each embedding.}
        \label{fig:enso_1.9}
\end{figure}

\newpage

\section{Additional computations: full persistence time-optimal PH cycle representatives}\label{A: full PH representatives}

\begin{figure}[h!]
\centering
\includegraphics[width=0.8\textwidth]{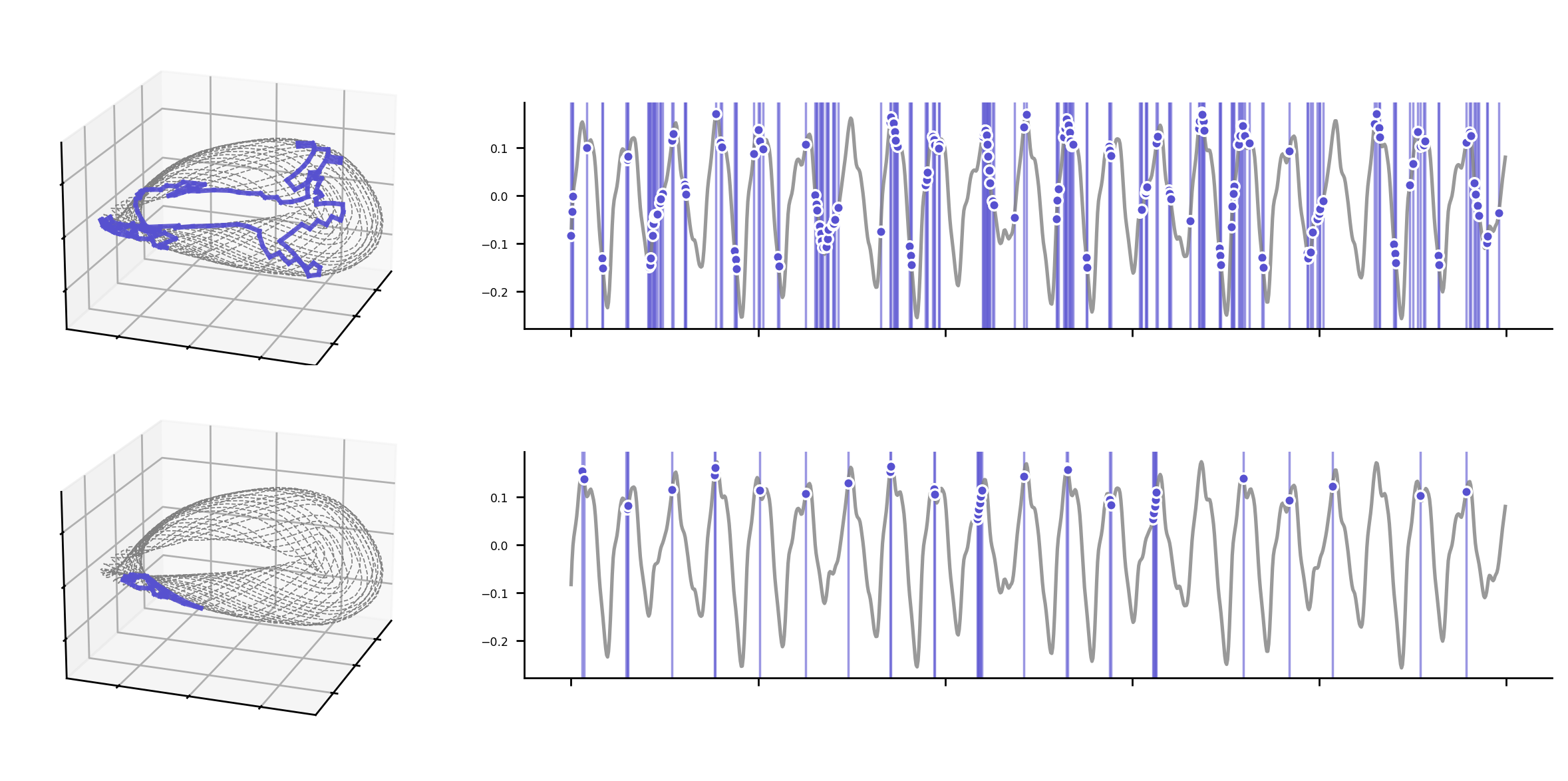}
\includegraphics[width=0.8\textwidth]{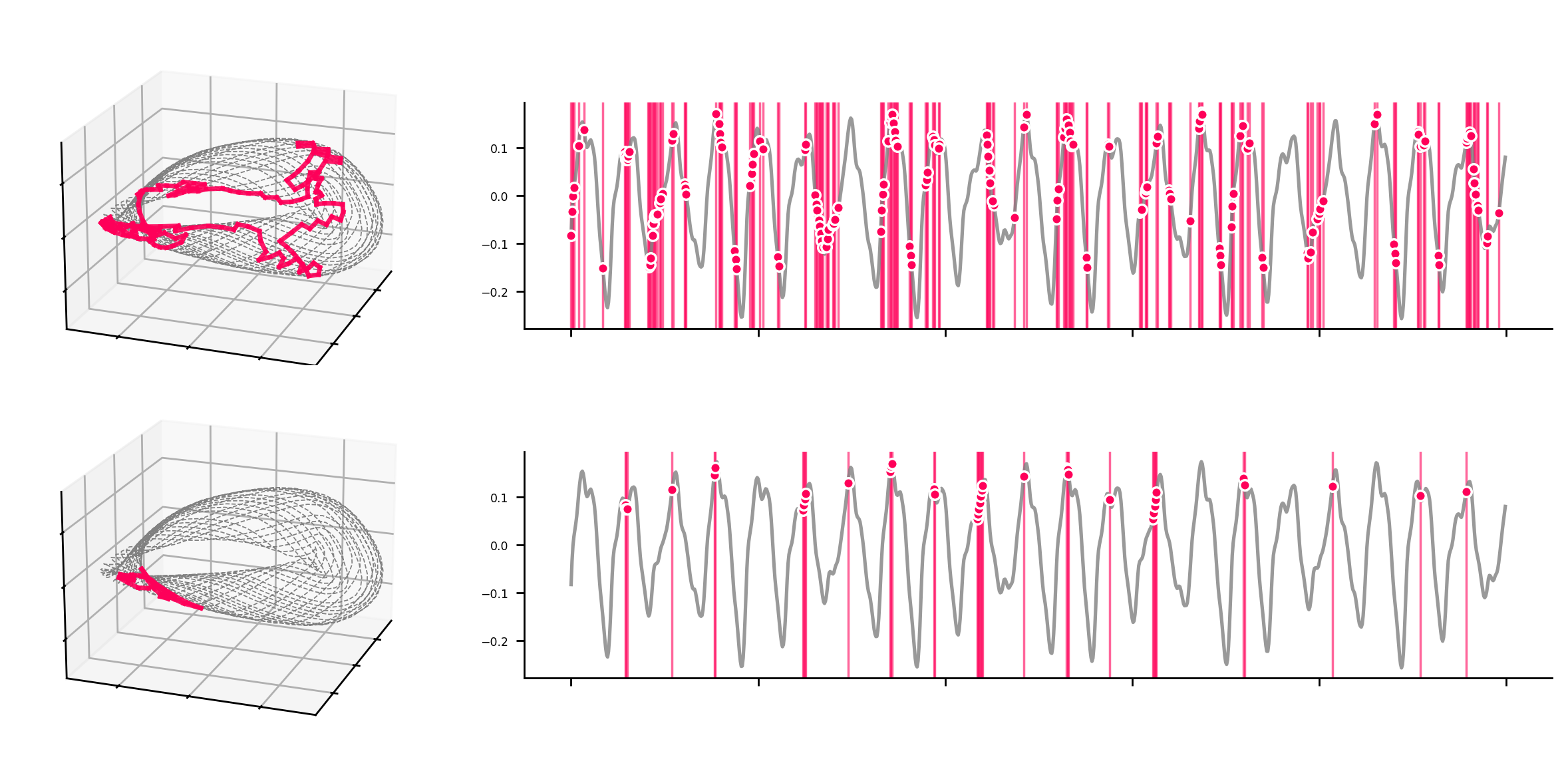}
\caption{Time representatives for the ENSO model with $\kappa=1.4$, simplex-based (top) and vertex-based (bottom), with full persistence}
\label{fig:example-enso14-representative-full}
\end{figure}

\begin{figure}[h!]
\centering
\includegraphics[width=0.8\textwidth]{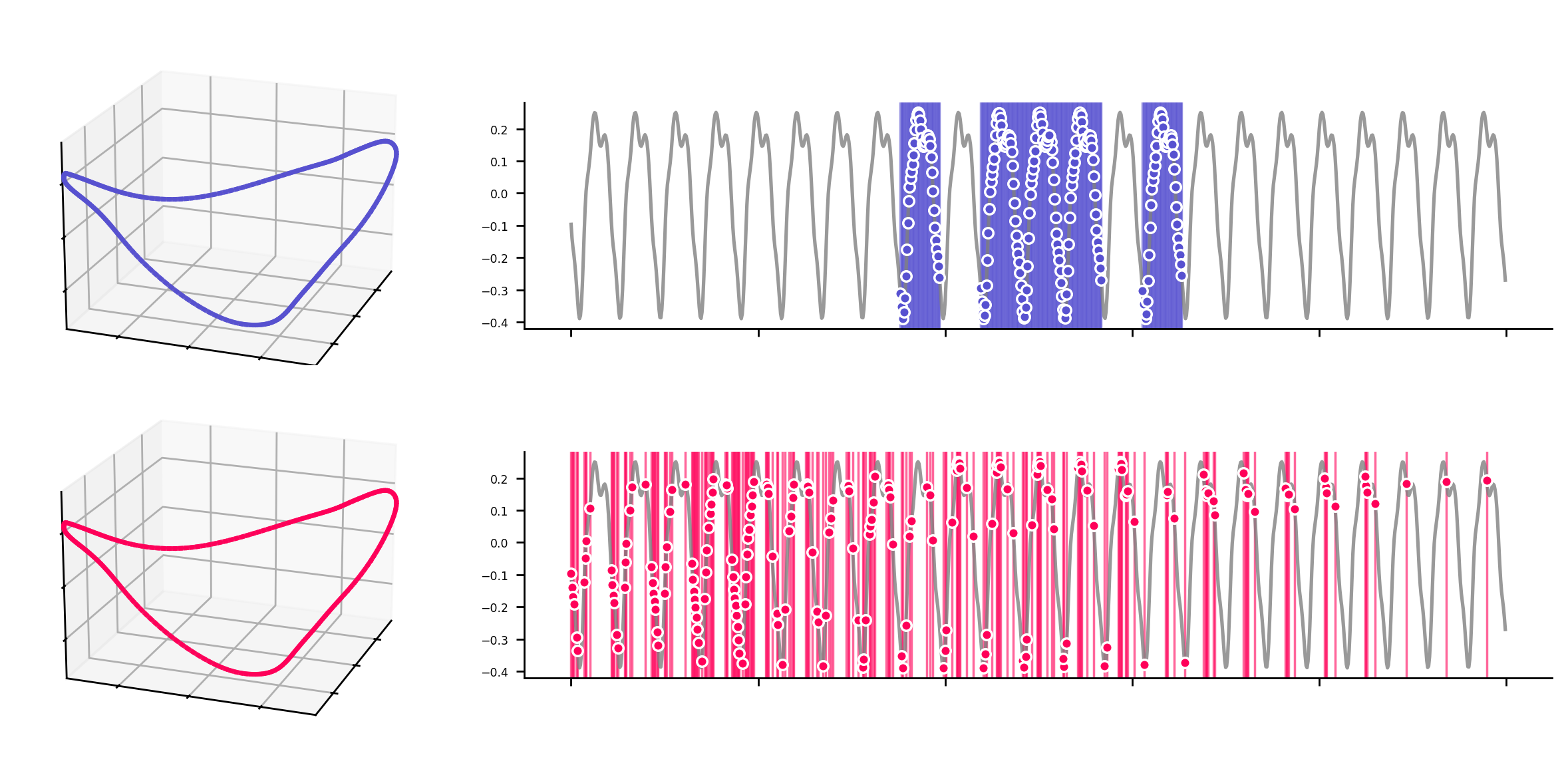}
\caption{Time representatives for the ENSO model with $\kappa=1.65$, simplex-based (top) and vertex-based (bottom), with full persistence}
\label{fig:example-enso165-representative-full}
\end{figure}

\begin{figure}[h!]
\centering
\includegraphics[width=0.8\textwidth]{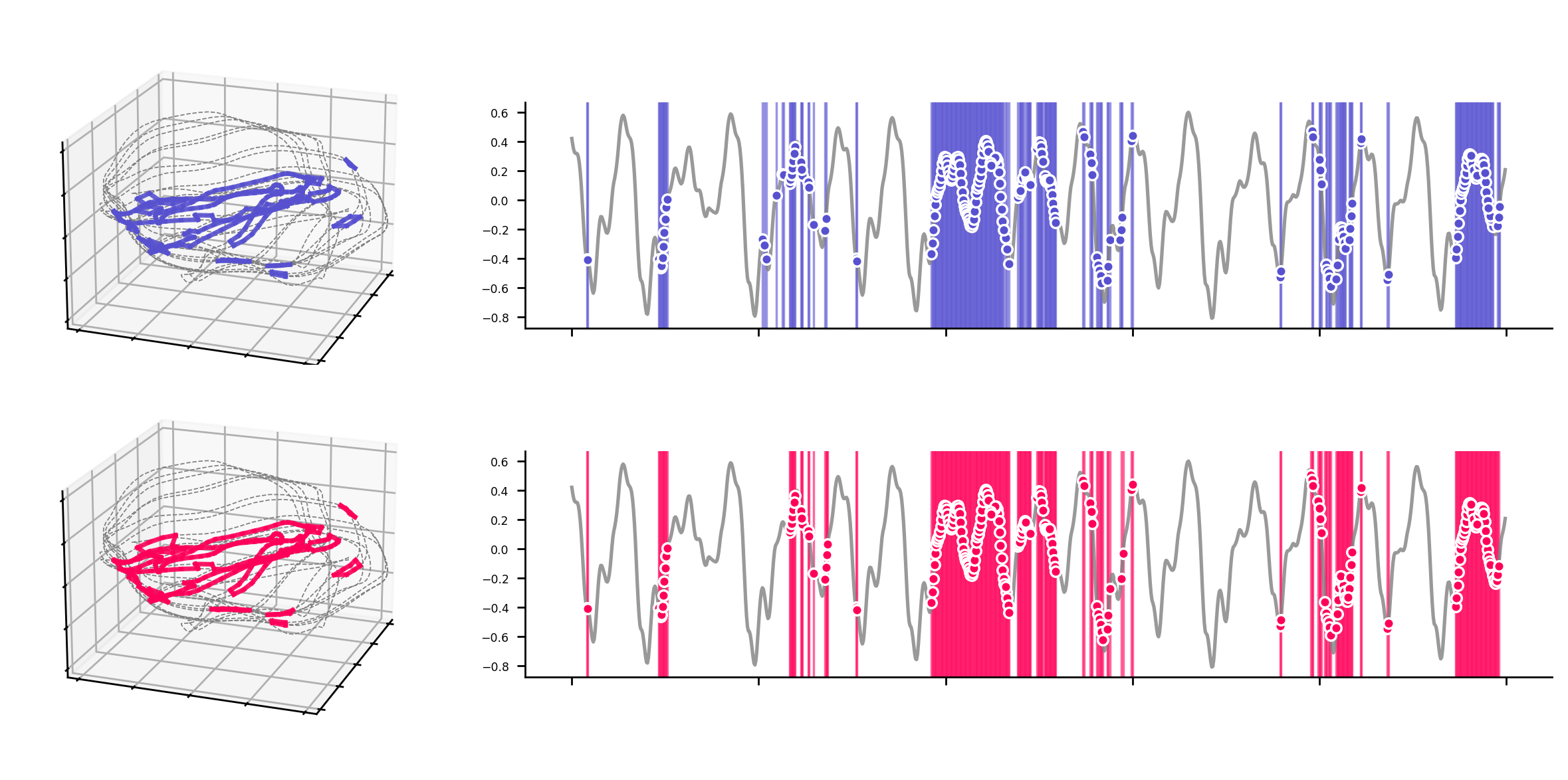}
\caption{Time representatives for the ENSO model with $\kappa=1.9$, simplex-based (top) and vertex-based (bottom), with full persistence.}
\label{fig:example-enso19-representative-full}
\end{figure}

\end{document}